\documentclass[12pt]{amsart}

\usepackage{framed,multirow}
\usepackage[mathlines]{lineno}
\usepackage[margin=1.2in]{geometry}

\usepackage{amssymb,amsthm}
\usepackage{latexsym}
\usepackage{dsfont}
\usepackage{hyperref}
\usepackage{amsfonts}
\usepackage{graphicx}
\usepackage{subfig}
\usepackage{multicol}
\usepackage{epstopdf}
\usepackage{overpic}
\usepackage{amsmath,esint}
\usepackage{epsfig}
\usepackage{tikz}
\usepackage{cleveref}
\usetikzlibrary{arrows}
\usetikzlibrary{fadings}
\usepackage{pgfplots}
\pgfplotsset{compat=1.10}
\usepgfplotslibrary{fillbetween}
\usetikzlibrary{patterns}
\usepackage{appendix}
\usepackage[utf8]{inputenc}
\usepackage{bbm}
\usepackage{enumitem}
\usepackage{mathrsfs}
\usepackage{verbatim}
\usepackage{graphicx}
\usepackage{algorithm}
\usepackage{tcolorbox}
\usepackage{algorithmicx}
\usepackage{algpseudocode}
\usepackage{mathtools}
\usepackage{pifont}
\usepackage{calc}
\usepackage{accents}
\usepackage{cancel}

\usepackage{ifpdf}
\usepackage{array}

\setenumerate{label=(\roman*),itemsep=3pt,topsep=3pt}

\newcommand{\RN}[1]{%
\textup{\uppercase\expandafter{\romannumeral#1}}%
}

\DeclarePairedDelimiter{\ceil}{\lceil}{\rceil}
\setenumerate{label=(\roman*),itemsep=3pt,topsep=3pt}

\newcommand{\cal}{\mathcal}

 \newcommand{\bbE}{\mathbb E}

\newcommand{\bbR}{\mathbb R}

\newcommand{\cC}{\mathcal C}

\newcommand{\cS}{\mathcal S}

\newcommand{\Rthree}{{\bbR^3}}
\newcommand{\Stwo}{{\cS^2}}

\newcommand{\rd}{\mathrm{d}}

\newcommand{\dt}{\Delta t}
\newcommand{\dx}{\Delta x}
\newcommand{\dv}{\Delta v}
\newcommand{\fij}{f_{i,j}}
\newcommand{\fijm}{f_{i-1,j}}
\newcommand{\fijp}{f_{i+1,j}}
\newcommand{\gij}{\lambda_{i,j}}
\newcommand{\gijm}{\lambda_{i-1,j}}
\newcommand{\gijp}{\lambda_{i+1,j}}
\newcommand{\np}{{m+1}}
\newcommand{\average}[1]{ \langle#1 \rangle}

\newtheorem{rem}{Remark}

\graphicspath{{Fig/}}

\begin{document}

\title{Adjoint Monte Carlo Method}

\author{Russel  Caflisch}
\email{caflisch@courant.nyu.edu}
\address{Courant Institute of Mathematical Sciences, New York University, New York, NY 10012.}

\author{Yunan Yang}
\email{yunan.yang@cornell.edu}
\address{Department of Mathematics, Cornell University, Ithaca, NY 14853.}

\begin{abstract} 
This survey explores the development of adjoint Monte Carlo methods for solving optimization problems governed by kinetic equations, a common challenge in areas such as plasma control and device design. These optimization problems are particularly demanding due to the high dimensionality of the phase space and the randomness in evaluating the objective functional, a consequence of using a forward Monte Carlo solver. To overcome these difficulties, a range of ``adjoint Monte Carlo methods'' have been devised. These methods skillfully combine Monte Carlo gradient estimators with PDE-constrained optimization, introducing innovative solutions tailored for kinetic applications. In this review, we begin by examining three primary strategies for Monte Carlo gradient estimation: the score function approach, the reparameterization trick, and the coupling method. We also delve into the adjoint-state method, an essential element in PDE-constrained optimization. Focusing on applications in the radiative transfer equation and the nonlinear Boltzmann equation, we provide a comprehensive guide on how to integrate Monte Carlo gradient techniques within both the optimize-then-discretize and the discretize-then-optimize frameworks from PDE-constrained optimization. This approach leads to the formulation of effective adjoint Monte Carlo methods, enabling efficient gradient estimation in complex, high-dimensional optimization problems.
\end{abstract}

\maketitle

\section{Introduction} \label{sec:Intro}
                                                                                                                                                                                                                                                                                                                                                                                                                                                                                                                                                                                                                                                                                                                                                                                                                                                                                                                                                                                                                                                                                                                                                                                                                                                                                                                                                                                                                                                                                                                                                                                                                                                                                                                                                                                                                                                                                                                                                                                                                                                                                                                                                                                                                                                                                                                                                                                                                                                                                                                                                                                                                                                                                                                                                                                                                                                                                                                                                                                                                                                                                                                                                                                                                                                                                                                                                                                                                                                                                                                                                                                                                                                                                                                                                                                                                                                                                                                                                                                                                                                                                                                                                                                                                                                                                                                                                                                                                                                                                                                                                                                                                                                                                                                                                                                                                                                                                                                                                                                                                                                                                                                                                                                                                                                                                                                                                                                                                                                                                                                                                                                                                                                                                                                                                                                                                                                                                                                                                                                                                                                                                                                                                                                                                                                                                                                                                                                                                                                                                                                                                                                                                                                                                                                                                                                                                                                                                                                                                                                                                                                                                                                                                                                                                                                                                                                                                                                                                                                                                                                                                                                                                                                                                                                                                                                                                                                                                                                                                                                                                                                                                                                                                                                                                                                                                                                                                                                                                                                                                                                                                                                                                                                                                                                                                                                                                                                                                                                                                                                                                                                                                                                                                                                                                                                                                                                                                                                                                                                                                                                                                                                                                                                                                                                                                                                                                                                                                                                                                                                                                                                                        The Monte Carlo method, a cornerstone in the field of computational mathematics, is a stochastic technique used for solving a wide array of mathematical problems which are often too complex for standard deterministic approaches~\cite{metropolis1949monte,rubinstein2016simulation}. 
The Monte Carlo method, as we know it today, was developed by John von Neumann and Stanis{\l}aw Ulam during World War II~\cite{von1947statistical,metropolis1949monte,eckhardt1987stan,metropolis1987beginning}, aiming to enhance decision-making in situations characterized by uncertainty. Its name is inspired by the famous casino town of Monaco, reflecting the method's fundamental reliance on the concept of randomness, akin to a roulette game. The Monte Carlo method is renowned for its versatility and efficacy in diverse domains ranging from physics to finance and engineering to artificial intelligence. At its core, the Monte Carlo method involves the use of random sampling to approximate solutions to problems that might be deterministic in principle but are practically infeasible to solve directly.

In many optimization problems arising from a wide array of applications where the Monte Carlo method is most effective, there is a strong need to compute the derivative of an objective function with respect to its parameters. This is crucial if one uses local optimization algorithms such as gradient descent and Newton's method to find the optimizer with fewer objective function evaluations~\cite{nocedal2006numerical}. A significant difficulty of directly computing the gradient is that the objective function is approximated by the Monte Carlo method, suffering from intrinsic randomness. Handling the randomness in the objective function gives rise to the field of Monte Carlo gradient estimation.

In computational mathematics, the technique of Monte Carlo gradient estimation emerges as a critical tool, leveraging stochastic processes to tackle complex optimization challenges~\cite{glasserman2013monte,mohamed2019monte,glynn1990likelihood,rubinstein1986score}. Building on the foundational Monte Carlo method, this approach broadens its application to gradient estimation. In particular, Monte Carlo gradient estimation excels in approximating the gradients of expectations, a task often impractical or impossible due to the complexity of the underlying expectations. This method skillfully overcomes such barriers by sampling from relevant distributions and using these samples for gradient estimation.

There are three primary strategies employed in Monte Carlo gradient estimators:
\begin{enumerate}
    \item the likelihood ratio method~\cite{glynn1990likelihood}, also known as the score function method~\cite{rubinstein1986score};
    \item  the pathwise derivative method~\cite{glasserman2013monte,glasserman1990gradient}, sometimes referred to as the reparameterization trick~\cite{kingma2013auto};
    \item the coupling method~\cite{pflug2012optimization,mohamed2019monte}.
\end{enumerate}
Each of these strategies possesses distinct characteristics and is suited to specific problem contexts. We will revisit these strategies comprehensively in Section~\ref{sec:background}, exploring their unique properties and applicability.

In this survey, our objective is to explore the application of Monte Carlo gradient estimation techniques in addressing optimization problems that emerge from kinetic theory. Notably, the Monte Carlo method has been a prominent choice for numerically solving kinetic equations, especially those with state variables existing in a seven-dimensional phase space~\cite{caflisch1980boltzmann}.  Alongside the development of mesh-based numerical methods, there has been significant interest in the Direct Simulation Monte Carlo (DSMC) approach and its various extensions~\cite{bird1970direct, nanbu1980direct, bobylev2000theory,pareschi2001introduction,pareschi2013interacting}. 

When dealing with design and control tasks modeled by kinetic equations, such as plasma control and devise design~\cite{einkemmer2023suppressing}, the challenge involves mathematically resolving a kinetic Partial Differential Equation (PDE)-constrained optimization problem. This classical approach typically requires adjusting parameters to align with reference data or to achieve specific desired properties, effectively minimizing an objective functional~\cite{hinze2008optimization}. The iterative update process in these problems involves the computation of gradients, requiring the solution of two PDEs: the original forward equation and its corresponding adjoint PDE, based on the adjoint-state method. In particular, the gradient is derived from either the optimize-then-discretize (OTD) or the discretize-then-optimize (DTO) frameworks (see Section~\ref{sec:PDE-constrained}). However, this may lead to prohibitive computational costs, particularly in high-dimensional PDEs.

To reduce the computational costs, one can consider using Monte Carlo solvers for gradient-based PDE-constrained optimizations. While there are established Monte Carlo methods for solving the complex forward kinetic PDE, such as the Boltzmann equation and the Vlasov--Poisson equation, a key research question arises: \newline
\begin{center}
\textit{How to effectively solve the corresponding adjoint equation using a Monte Carlo method? }
\end{center}

The primary challenge in developing a Monte Carlo numerical solver for the adjoint equation lies in reconciling the Monte Carlo approach with the computation of gradients. In essence, the Monte Carlo forward solver computes the weak form solution of a given PDE, which is  a weighted linear combination of Dirac delta measures. However, the gradient calculations for PDE-constrained optimization problems typically involve multiplying two PDE solutions – one forward and one adjoint. Utilizing the same or a similar Monte Carlo solver for the adjoint PDE as used for the forward PDE presents a problem; the gradient requires the multiplication of two sets of Dirac delta measures over the same phase space, which is not well-defined since Dirac delta is not technically a function. Therefore, it is crucial to carefully design an appropriate adjoint Monte Carlo method, which also requires a thorough understanding of the nature of adjoint variables. This survey focuses on this exact issue.

\begin{figure}
\begin{center}
\begin{tikzpicture}
    \draw[black,thick] (0,0) circle (0.3*\textwidth);
    \draw[black,thick] (\textwidth/3,0) circle (0.3*\textwidth);
    \node at (0,4) {\textbf{MC Gradient}};
    \node at (0,3.5) {\textbf{Estimation}};
    \node at (-2,1.5) {1.~Score Function};
    \node at (-2.0, 1) {Method (see Sec.~\ref{subsubsec:score})};
    \node at (-1.6,0) {2.~Pathwise Derivative};
    \node at (-2,  -0.5) {Method (see Sec.~\ref{subsubsec:pw_grad})};
    \node at (-1.8,-1.5) {3.~Coupling Method};
    \node at (-2, -2) {(see Sec.~\ref{subsubsec:coupling})};
    \node at (5.25,4) {\textbf{PDE-Constrained}};
    \node at (5.25,3.5) {\textbf{Optimization}};
    \node at (6.5,1) {1.~OTD vs.~DTO};
    \node at (6.5,0.5) {(see Sec.~\ref{subsec:OTD_DTO})};
    \node at (7.2,-0.5) {2.~Adjoint-State Method};
    \node at (6.5,-1) {(see Sec.~\ref{subsec:adjoint})};
    \node at (2.6,2.2) {\textbf{Adjoint MC}};
    \node at (2.6,1.7) {\textbf{Method}};
    \node at (2.6,0.5) {Case Study I:};
    \node at (2.6,0) {Radiative Transport};
    \node at (2.6,-0.5) {Eqn.~(see Sec.~\ref{sec:rte})};
    \node at (2.6,-1.5) {Case Study II:};
    \node at (2.6,-2) {Boltzmann Eqn.~}; 
    \node at (2.6,-2.5) {(see Sec.~\ref{sec:Boltzmann})};
\end{tikzpicture}
\end{center}
\caption{The diagram illustrates the main structure of this survey paper. The adjoint Monte Carlo (MC) method integrates effective techniques from Monte Carlo gradient estimation and problem-solving frameworks in PDE-constrained optimization to tackle challenging optimization problems based on various kinetic PDEs.\label{fig:overview}}
\end{figure}
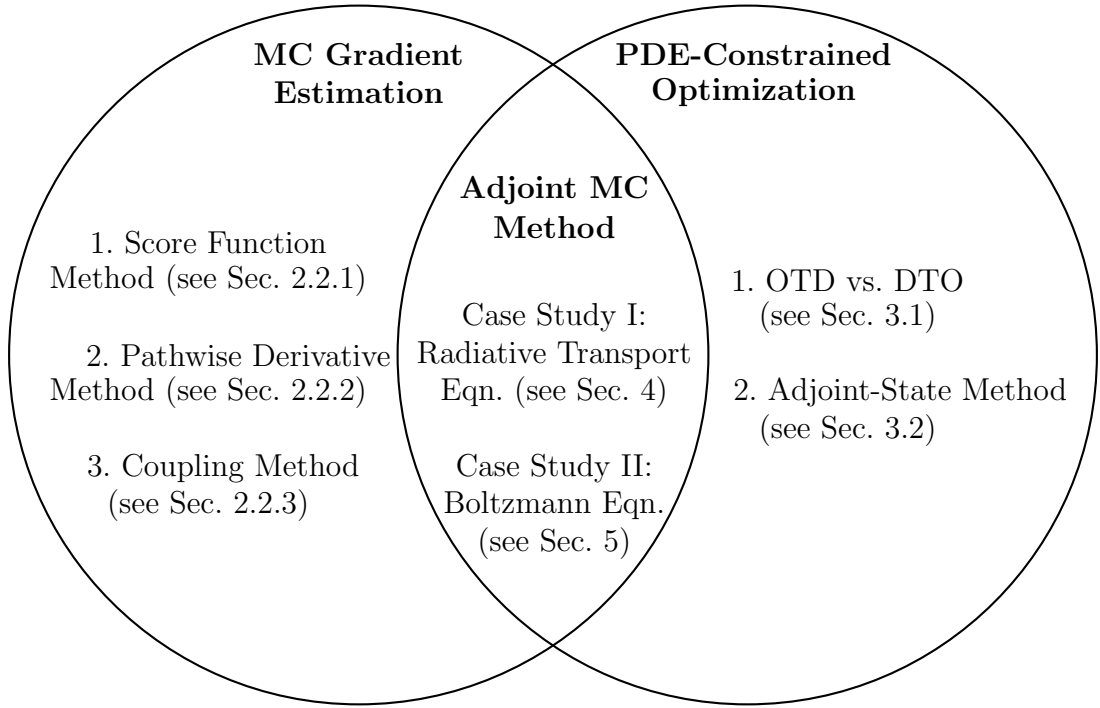

In this survey, we present a  general framework for the adjoint Monte Carlo method, illustrating its application in two case studies: the Radiative Transfer Equation (RTE)~\cite{li2022monte} and the Boltzmann Equation~\cite{caflisch2021adjoint,yang2023adjoint}. The method skillfully combines Monte Carlo gradient estimation techniques, prevalent in mathematical finance and machine learning~\cite{glasserman2013monte,mohamed2019monte}, with established solutions for PDE-constrained optimization~\cite{herty2007optimal}. It introduces a novel approach to handling randomness in the discretization and approximation of forward equations, which is beyond the typical study of PDE-constrained optimization, necessitating new theoretical and computational strategies for the adjoint PDE. An overview diagram of this method is also provided in~\Cref{fig:overview}.

We aim to establish a framework for formulating an adjoint Monte Carlo method, showcased through two detailed case studies: the radiative transfer equation (RTE)~\cite{li2022monte} and the Boltzmann equation~\cite{caflisch2021adjoint,yang2023adjoint}. This exploration will not only address the theoretical aspects but also provide practical insights into its implementation. In a nutshell, the adjoint Monte Carlo method designed for various kinetic PDEs skillfully integrates mature techniques from Monte Carlo gradient estimation, which are widely used in mathematical finance and machine learning~\cite{glasserman2013monte,mohamed2019monte}, with solution frameworks for PDE-constrained optimization problems. Unlike the typical solution procedure in PDE-constrained optimization, the discretization and approximation of the forward (state) equation in the context considered here involve randomness, which requires innovations in both theoretical understandings and computational algorithms for the adjoint PDE. An overview diagram for the adjoint Monte Carlo method is presented in~\Cref{fig:overview}.

Since its first proposal in 2021~\cite{caflisch2021adjoint}, the adjoint Monte Carlo method~\cite{yang2023adjoint,li2022monte} has seen impactful use, notably in topology optimization for rarefied gas flow channel design~\cite{guan2023topology,yuan2023efficient}, a spatially inhomogeneous Boltzmann-equation constrained optimization problem. A key challenge of the adjoint Monte Carlo method is its high memory requirement, which led to the development of a reversible extension for permuted congruential pseudorandom number generators to reduce memory usage~\cite{lovbak2023reversible}. This further influences the work of adjoint sensitivity calculation in plasma edge codes~\cite{carli2023algorithmic}. The adjoint Monte Carlo method for kinetic applications has gathered much research interest in recent years, which warrants the timeliness of this survey. 

The rest of the paper is organized as follows. We first introduce essential mathematical and statistical background in Section~\ref{sec:background}, where we will present an overview of the field of the Monte Carlo gradient and its mainstream strategies. In Section~\ref{sec:PDE-constrained}, we switch gears to revisit PDE-constrained optimization problems and the general mathematical frameworks to tackle them. In particular, we re-derive the adjoint-state method, a key ingredient to guarantee the scalability of PDE-constrained optimization for real-world applications. The central theme of this survey paper is to showcase how one can combine the Monte Carlo gradient estimation methods in Section~\ref{sec:background} with strategies that solve the PDE-constrained optimization problem in Section~\ref{sec:PDE-constrained}, two separate research topics, together, and innovate a new class of method, termed as the ``Adjoint Monte Carlo Method'', for solving optimization problems regarding kinetic PDEs. In Section~\ref{sec:rte}, we give an example of the radiative transport equation (RTE), a linear mean-field equation. To compute the derivative of the objective functional constrained by the RTE, we show how to derive and utilize the adjoint Monte Carlo method based on the OTD and DTO approaches, respectively. In Section~\ref{sec:Boltzmann}, we look closely at the Boltzmann equation, a nonlinear integro-differential equation but not a mean-field equation. The binary collision in the Boltzmann equation requires careful treatment, and we obtain the adjoint Monte Carlo method in the DTO approach. With the two examples, we have outlined the steps to construct an adjoint Monte Carlo scheme. Conclusions and more discussions follow in Section~\ref{sec:conclusion}.

\section{Monte Carlo Gradient Estimation}\label{sec:background}
The term ``Monte Carlo gradient'' refers to the estimation of the gradient of a given objective function using Monte Carlo methods. Consider a general objective function in the form of
\begin{equation}\label{eq:main_obj}
    J(\theta) = \int f(x) p(x,\theta) \rd x = \bbE_{X\sim p(x,\theta)}\left[f(X) \right],
\end{equation}
where $f(x)$ is an observable function, assumed to be independent of $\theta$, and $p(x,\theta)$ is the input density function corresponding to a given distribution that depends on the parameter $\theta$. The optimization problem aims to find an optimal parameter $\theta$ that minimizes the objective function $J(\theta)$. If one wants to use optimization algorithms such as steepest descent, conjugate gradient descent, and Newton's method~\cite{nocedal2006numerical}, it is essential to have a computationally feasible way of obtaining the gradient of $J(\theta)$. 

Often, we do not have direct access to the analytic form of the density function $p(x,\theta)$, or it is challenging to use other deterministic quadrature rules to approximate the integral in~\eqref{eq:main_obj}. In such cases, the objective function~\eqref{eq:main_obj} can be approximated through the Monte Carlo integration:
\begin{equation}\label{eq:MC_obj}
    J(\theta) \approx   \widehat{J}(\theta) = \frac{\rho}{N} \sum_{i=1}^N f(X_i),\quad X_1,\ldots, X_N \sim p(x,\theta),
\end{equation}
where $X_1,\ldots, X_N $  are independent and identically distributed (i.i.d.) samples, $\rho = \int p(x,\theta) \rd x$ is the total mass and is often independent with respect to the parameter $\theta$. For probability distributions, we always have $\rho = 1$. We remark that $\widehat{J}(\theta)$ given in~\eqref{eq:MC_obj} is an unbiased estimator as it satisfies 
$$\bbE\left[ \widehat{J}(\theta) \right] = J(\theta),$$
where the expectation is taken with respect to all the i.i.d.~samples $X_1,\ldots, X_N$. 

Based on the Central Limit Theorem, we have~\cite{caflisch1998monte}
\[
J(\theta)  - \widehat{J}(\theta)  \approx \sigma(\theta) N^{-1/2} \, \mathbf{r}\,,
\]
where the random variable $\mathbf{r}\sim N(0,1)$, the 1D standard normal distribution, and $$\sigma(\theta) = \sqrt{\int \left(f - J(\theta) \right)^2 \, p(x,\theta) \rd x}$$ 
is the standard deviation, depending on the parameter $\theta$.

In this survey paper, the problem of interest is not to estimate $J(\theta) $ but instead to evaluate (or approximate) the gradient of $J(\theta)$ with respect to the parameter $\theta$, 
\begin{equation}\label{eq:true_grad}
G(\theta) := \nabla_\theta J(\theta)\,.
\end{equation}
There are two main scenarios where the Monte Carlo gradient estimation is useful. First, in some contexts, especially in high-dimensional spaces or complex models, it is challenging or computationally expensive to directly compute the gradient $G(\theta)$. Monte Carlo gradient estimation provides a stochastic approach to approximate these gradients. Second, in many cases, the objective function evaluation can be stochastic in nature due to either the computational complexity of directly evaluating $J(\theta)$  or the lack of access to the full distribution function $p(x,\theta)$. In such scenarios, one has to approximate the objective function $J(\theta)$  through certain Monte Carlo methods. As a result, one only has access to the random variable $\widehat{J}(\theta)$ instead of deterministic value $J(\theta)$. To examine the sensitivity of the approximated random objective function value $\widehat{J}(\theta)$ to the parameter $\theta$, it is then equivalent to examining the sensitivity of the samples 
$$
X_1,X_2,\ldots, X_N \sim p(x;\theta)
$$
in the Monte Carlo method with respect to the parameter $\theta$ in the underlying distribution. We will use $\widehat{G}(\theta)$ to denote the Monte Carlo gradient estimation of $G(\theta)$.

\subsection{Overview}

Monte Carlo gradient estimation, sometimes referred to as the ``differentiable Monte Carlo'', is a topic of interest and an active field of research in a wide range of subjects~\cite{mohamed2019monte}. In reinforcement learning~\cite{kaelbling1996reinforcement}, one research topic is to determine the optimal policy that has a maximum reward. The objective is often an expected reward where the reward is a random variable, similar to the form of $J(\theta)$ in~\eqref{eq:main_obj}. The so-called policy gradient, the gradient of this expectation with respect to the policy parameters in reinforcement learning, is often estimated using Monte Carlo methods. The well-known REINFORCE algorithm is a classic example~\cite{williams1992simple}, which utilized the score-function method, which we will discuss in detail later in Section~\ref{subsubsec:score}.

In variational Bayesian inference, the evidence lower bound (ELBO) is a fundamental concept, which is a method used to approximate complex posterior distributions. The goal of variational inference is to maximize the ELBO with respect to the parameters of the variational distribution. This is equivalent to minimizing the divergence between the true posterior and the variational distribution. The variational problem is 
often optimized using Monte Carlo gradient estimators~\cite{kingma2013auto}, especially when the model involves latent variables. 
The Monte Carlo gradient estimation provides an unbiased estimate of the true gradient. This makes it a powerful tool for optimizing the ELBO, especially when the exact gradient computation is intractable.

Monte Carlo methods are widely used in mathematical finance for valuing complex financial derivatives and risk management. Some financial derivatives have payoffs that are challenging to value using traditional closed-form solutions or numerical methods like finite difference methods.
Monte Carlo simulation can be used to simulate many possible sample paths of the underlying asset(s) and then average the payoffs across these paths to estimate the derivative's value.
The gradient (sensitivity) of this value with respect to model parameters can be estimated using the Monte Carlo gradient, e.g., the ``Greeks'' of financial instruments with respect to underlying parameters~\cite[Sec.~7.2]{glasserman2013monte}.

Monte Carlo methods offer a notable advantage in their natural capacity for parallelization. This characteristic stems from the methods' reliance on numerous i.i.d.~random samples to approximate solutions. Since these sampling activities are independent of each other, they can be effectively distributed across various processors or computing nodes. This ability to execute in parallel not only speeds up the computational process but also enhances the scalability and efficiency of Monte Carlo methods, making them particularly suitable for tackling large-scale problems in high-performance computing environments.

A significant challenge in utilizing Monte Carlo methods for optimization tasks, as previously noted, lies in calculating the gradients of expectations estimated through these methods, inherently random in nature. Accurately determining these gradients is essential for optimization in fields such as reinforcement learning, Bayesian inference, and simulation-based design. The stochastic nature of Monte Carlo gradient estimation means that the estimated gradient can have high variance. Techniques like quasi Monte Carlo method~\cite{caflisch1998monte}, importance sampling~\cite{tokdar2010importance}, antithetic variates~\cite{kroese2013handbook}, or control variates~\cite{rubinstein1985efficiency} can be used to reduce this variance. Depending on the sampling method, there might be bias in the gradient estimates, i.e., $\bbE[\widehat{G}(\theta)]\neq G(\theta)$. It is crucial to ensure that the estimator is unbiased or that the bias is controlled and understood. This topic has attracted much research endeavor~\cite{mark2001bias,chen2005stochastic,tucker2017rebar,nemeth2021stochastic}.

\subsection{Methods}
In the following subsections, we will introduce three primary methods of Monte Carlo gradient estimation, each widely employed across diverse fields.

\subsubsection{Likelihood Ratio Method / Score Function Method}\label{subsubsec:score}

The Monte Carlo gradient estimator based on the score function method stands as one of the most versatile gradient estimators. It has a wide range of applications. It is known as the score function estimator in~\cite{rubinstein1986score,kleijnen1996optimization}, and referred to as the likelihood ratio method in~\cite{glynn1990likelihood}. The method is also called the REINFORCE estimator in reinforcement learning~\cite{williams1992simple}. Next, we will briefly introduce its derivation and comment on the properties.

Consider the parameterized density function $p(x,\theta)$. The score function is the derivative of the logarithm of $p(x,\theta)$, i.e.,
\begin{equation}\label{eq:score_def}
    \nabla_\theta  \log p(x,\theta)\,,
\end{equation}
and it satisfies the relationship
\begin{equation}\label{eq:score}
\nabla_\theta p(x,\theta) = p(x,\theta) \nabla_\theta  \log p(x,\theta) \,.
\end{equation}
If $X$ is a random variable with its distribution characterized by the density function $p(x,\theta)$, we have
another important property regarding the score function:
\begin{equation}\label{eq:score_0}
    \bbE_{X\sim p(x,\theta)} \left[ \nabla_\theta  \log p(X,\theta) \right] = \int \frac{\nabla_\theta p(x,\theta) }{p(x,\theta)} p(x,\theta) \rd x = \nabla_\theta \int  p(x,\theta) \rd x = \nabla_\theta 1 = 0.
\end{equation}
That is, as a random variable, $ \nabla_\theta  \log p(X,\theta)$ has mean zero. If $\theta \in \bbR^m$, $\nabla_\theta  \log p(X,\theta) $ is a vector-valued random variable.
We proceed to compute its covariance matrix $C(\theta)$. The $ij$-th entry of $C(\theta)$ is
\begin{eqnarray}
   C_{ij}(\theta) &=& \int  \frac{\partial \log p(x,\theta) }{\partial \theta_i}   \,\frac{\partial \log p(x,\theta) }{\partial \theta_j}   p(x,\theta)\rd x \nonumber   \\
   &=&   - \bbE_{X\sim p(x,\theta) }\left[ \frac{\partial^2}{\partial \theta_i \partial\theta_j} \log   p(X,\theta)  \right]\,.
   \label{eq:score_var} 
\end{eqnarray}
The last term is the so-called Fisher information. The second equality holds because $$\nabla_\theta \left(  \bbE_{X\sim p(x,\theta)} \left[ \nabla_\theta  \log p(X,\theta) \right]  \right) = 0\,,$$
as a result of~\eqref{eq:score_0}.

The score function plays a crucial role in statistical estimation and machine learning, serving as a foundational tool for parameter estimation and model evaluation. It is essentially the gradient of the log-likelihood function with respect to the parameters, providing critical information about how the likelihood changes with respect to these parameters. The importance of the score function lies in its ability to identify the most probable parameters that best fit the data, thereby maximizing the likelihood. In complex models, the score function guides optimization algorithms toward the best parameter values, significantly impacting the accuracy and efficiency of model fitting. Additionally, it is instrumental in constructing various statistical tests and confidence intervals, making it a key component in inferential statistics and model diagnostics.

Given the properties of the score function, we next describe how to use it to estimate the gradient defined in~\eqref{eq:true_grad}. We first reformulate $G(\theta)$ as follows.
\begin{eqnarray}
    G(\theta) &=&  \nabla_\theta \left( \int f(x) p (x,\theta) \rd x \right) \nonumber \\
    &=& \int f(x) \nabla_\theta p (x,\theta) \rd x \nonumber \\
    &=& \int f(x) (\nabla_\theta \log p (x,\theta)) \, p(x,\theta) \rd x  \nonumber \\
    &=& \bbE_{X\sim p (x,\theta)} \left[  f(X) \nabla_\theta \log p (X,\theta)\right]\,,\label{eq:score_grad1}
\end{eqnarray}
where we used the identity in~\eqref{eq:score}. As a result,  to generate an estimation of $G(\theta)$, denoted by $\widehat {G}$, we can use Monte Carlo integration to approximate the last term in~\eqref{eq:score_grad1}:
\begin{equation}\label{eq:score_grad2}
    \widehat {G} \approx \frac{1}{N} \sum_{i=1}^N f(X_i) \, \nabla_\theta \log p (X_i,\theta)\,,
\end{equation}
where $X_1,\ldots, X_N$ are i.i.d.~samples following the distribution $p (x,\theta)$. Here, we assume that $\rho  = \int p (x,\theta) \rd x = 1$. 
The formula outlined in~\eqref{eq:score_grad2} is the score method gradient estimator, one of the most classic techniques for gradient estimation. We can also verify that the estimator $\widehat{G}$ is unbiased:
\[
\bbE_{X_i\sim p(x,\theta),\text{i.i.d.}}[\widehat {G} ] = G(\theta)\,.
\]

Recall that in the beginning of~\Cref{sec:background}, we remarked that estimating the Monte Carlo gradient is equivalent to examining the sensitivity of the samples $\{X_i\}$ on the parameter $\theta$. The ``Score Function Method'' achieves this by first evaluating the score function at those sample positions, i.e., $\nabla_\theta \log p (X_i,\theta)$, and then multiplying the contribution by a weight $f(X_i)$. The final weighted average of the score function values yields the true parameter sensitivity.

Note that we interchanged the order between differentiation and integration in the reformulation of $G(\theta)$ in~\eqref{eq:score_grad1}. To justify this step, we need extra conditions on the density function $p(x,\theta)$ and the integrand $f(x)$~\cite{l1995note,kleijnen1996optimization}. One necessary condition is that $p(x,\theta)$ is continuously differentiable with respect to $\theta$ and is also absolutely
continuous with respect to $\theta$ at the boundary of the support~\cite{pflug2012optimization,glasserman2013monte}.

\subsubsection{Pathwise Derivative Method / Reparameterization Trick}\label{subsubsec:pw_grad}
The pathwise derivative method~\cite{glasserman2013monte} is also sometimes referred to as the process derivative~\cite{pflug2012optimization,glasserman1990gradient}, stochastic backpropagation~\cite{rezende2014stochastic}, and the reparameterization trick~\cite{kingma2013auto}. The pathwise derivative method is often used to make the gradient estimation more stable. This method involves reparameterizing the latent variables in a way that allows the gradient to be taken with respect to the parameters of the variational distribution directly rather than through the sampling process.

To begin with, we consider a specific way of sampling the target probability measure  $\mu_\theta$ with density $p(x,\theta)$ where $\rd\mu(x) = p(x,\theta) \rd x$. Suppose there exists  a parameter $\theta$-independent distribution $\pi$ such that
\begin{equation}\label{eq:push-forward}
   \mu_\theta = T_\theta \sharp  \pi,\quad T_\theta: \bbR^n \mapsto \bbR^n\,.
\end{equation}
Here, the function $T_\theta$ is often referred to as the push-forward map. It turns one measure $\pi$ to another measure $\mu$ and is also mass-preserving: $$\mu_\theta(A) =\pi \left(T_\theta^{-1}(A) \right) \, $$ 
where $A$ is an arbitrary Borel measurable set and $T_\theta^{-1}(A)$ denotes its preimage.

If the relationship~\eqref{eq:push-forward} is satisfied, we can sample from $\mu_\theta$ by sampling $\pi$. To do so, we first obtain $N$ i.i.d.~samples from $\pi$, denoted by $Y_1,\dots, Y_N$, and transform these samples through the deterministic function $T_\theta$:
\begin{equation}\label{eq:path_sample}
    X_i = T_\theta(Y_i),\quad i = 1,\ldots, N\,.
\end{equation}
Then the resulting samples $X_1,\ldots, X_i,\ldots, X_N \sim \mu_\theta$. If $\pi$ is also absolutely continuous with respect to the Lebesgue measure in $\bbR^n$ with the density function $q(x)$, and moreover, $T_\theta$ is a bijection, then we have the change of variable formula:
\[
q(x) = p\left( T_\theta (x) ,\theta \right) |\nabla_x T_\theta|\,,
\]
where $|\nabla_x T_\theta|$ represents the determinant of the Jacobian matrix $\nabla_x T_\theta$. 

The parameterized push-forward map $T_\theta$ is also considered as a path function. For example, when sampling a Gaussian random variable following $N(3,\theta^2)$ with $\theta > 0$, we can first sample from the standard normal distribution $N(0,1)$, followed by proper translation and dilation through the linear map $T_\theta (x) = \theta \, x + 3 $.

A direct consequence of this pathwise sampling is the Law of the Unconscious Statisticians (LOTUS)~\cite{federer2014geometric}:
\begin{equation}\label{eq:LOTUS}
    \bbE_{Y\sim \pi} \left[ f\left(T_\theta(Y) \right) \right] = \bbE_{X\sim \mu_\theta} \left[ f\left(X \right) \right]\,.
\end{equation}
LOTUS tells us that we can perform (Monte Carlo) integration without
knowing the underlying probability distribution as long as we are given the corresponding sampling path $T_\theta$ and the base distribution $\pi$. Using this property, we obtain another reformulation of the  gradient $G$ that is different from~\eqref{eq:score_grad1}:
\begin{eqnarray}
    G &=& \nabla_\theta \left(\bbE_{X\sim \mu_\theta} \left[ f\left(X \right) \right]\right)  \nonumber \\
    &=& \nabla_\theta \left( \bbE_{Y\sim \pi} \left[ f\left(T_\theta(Y) \right) \right] \right) \nonumber \\
    &=&  \bbE_{Y\sim \pi} \left[ \nabla_\theta f\left(T_\theta(Y) \right) \right]\nonumber \\
    &=& \bbE_{Y\sim \pi} \left[   \nabla_\theta T_\theta(Y) ^\top \nabla_x f(x)\big|_{x=T_\theta(Y)} \right]  \label{eq:pw_grad1} \,.
\end{eqnarray}
The gradient estimator is then a Monte Carlo approximation of~\eqref{eq:pw_grad1}, given by
\begin{equation}\label{eq:pw_grad2}
    \widehat{G}  = \frac{1}{N}\sum_{i=1}^N \nabla_\theta T_\theta\big(Y_i) ^\top  \, \nabla_x f\left(T_{\theta}(Y_i)\right),\quad Y_i \sim \pi\,,\,\, i = 1,\ldots, N\,. 
\end{equation}

\begin{figure}
\begin{center}
\usetikzlibrary{arrows.meta}
\begin{tikzpicture}
    \draw [very thick] (-2,0) circle (1.8);
    \node at (-2,1) {$T_\theta(\boldsymbol{\cdot})$};
    \draw [gray, very thick, dashed] (-2,-0.75) 
 circle (0.5cm);
    \node[red] at (-2,-0.75) {$\theta$};
    \draw [-stealth](-2,0.75) -- (-2,-0.2);
    \node at (-1.2,0.3) {depends};
    \node at (-1.2,-0.1) {on};
    \draw [-stealth](-2,0.75) -- (-2,-0.2);
    \node(Y1) at (-5.5,2) {$Y_1$};
    \node(Y2) at (-5.5,1) {$Y_2$};
    \node(YY) at (-5.5,0) {$\vdots$};
    \node(YN) at (-5.5,-1.5) {$Y_N$};
    \draw [-latex] (Y1) -- (-3.6,0.9);
    \draw [-latex] (Y2) -- (-3.8,0.2);
    \draw [-latex] (YN) -- (-3.8,-0.5);
    \node(X1) at (1.5,2) {$X_1$};
    \node(X2) at (1.5,1) {$X_2$};
    \node(XX) at (1.5,0) {$\vdots$};
    \node(XN) at (1.5,-1.5) {$X_N$};
    \draw [-latex] (-0.4,0.9) -- (X1);
    \draw [-latex] (-0.2,0.2)-- (X2);
    \draw [-latex] (-0.2,-0.5)-- (XN);
    \node(fX1) at (3,2) {$f(X_1)$};
    \node(fX2) at (3,1) {$f(X_2)$};
    \node(fXX) at (3,0) {$\vdots$};
    \node(fXN) at (3,-1.5) {$f(X_N)$};
    \draw [-latex] (X1) -- (fX1);
    \draw [-latex] (X2)-- (fX2);
    \draw [-latex] (XN) -- (fXN);
    \node(J) at (7,0) {$ \widehat{J} = \frac{1}{N} \sum_{i=1}^N f(X_i)$};
    \draw[blue, very thick] (5,-0.5) rectangle (9,0.5);
    \draw [-latex] (fX1) -- (4.9,0.45);
    \draw [-latex] (fX2) -- (4.9,0.25);
    \draw [-latex] (fXN) -- (4.9,-0.35);
\end{tikzpicture}
\end{center}
\caption{A diagram showing the dependence of the approximated objective function $\widehat{J}$  on the samples $\{X_i\} \sim \mu_\theta$, which are determined by the parameter $\theta$ through the push-forward map $T_\theta(\boldsymbol{\cdot})$ \textbf{path by path} in the ``Pathwise Derivative Method'' (see Section~\ref{subsubsec:pw_grad}).\label{fig:pwd_diagram}}
\end{figure}
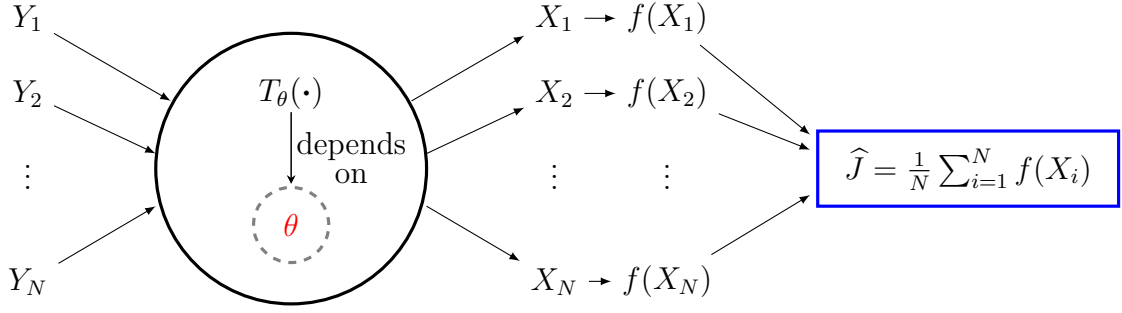

The pathwise gradient estimation is applicable when the function inside the expectation is differentiable with respect to its parameter $\theta$. It involves directly differentiating the simulation, and the formula inherits the structure of a chain rule:
\begin{center}
From $\theta$ to $T_\theta(X)$ (i.e., $\nabla_\theta T_\theta$), and from $T_\theta(X)$ to $f(T_\theta(X))$ (i.e., $\nabla_x f\left(x\right)|_{x=T_\theta(X)}$). \\
\end{center}
Different from the score function method discussed in Section~\ref{subsubsec:score}, the contribution of the gradient is traced path by path as shown in the diagram in~\Cref{fig:pwd_diagram}. The dependence of the samples $\{X_i\}_{i=1}^N$ on the parameter $\theta$ is evident through the deterministic path function (i.e., the push-forward map) $T_\theta$. This is the fundamental reason why the simple chain rule is applicable here.

The pathwise method is, however, limited to distributions that can be reformulated as~\eqref{eq:push-forward}, the push-forward of a baseline parameter-independent distribution. Moreover, the push-forward map has to be differentiable in terms of the parameters. Often, the process of sampling $X$ from $p(x,\theta)$ and the process of computing the gradient of the samples with respect to the parameter $\theta$, i.e., the computation of the $\nabla_\theta X$, are coupled together. It is possible to decouple these two processes by sampling using a favorite and efficient method, and then approximate $\nabla_\theta X$ in another way~\cite[Sec.~5.3]{mohamed2019monte}. Similar to the derivation of the score function estimator, we also used the interchange between differentiation and integration in~\eqref{eq:pw_grad1}. A sufficient condition requires $f(x)$ to be differentiable. If the interchange is justified, we can again see that the pathwise gradient estimator $\widehat{G}$ given in~\eqref{eq:pw_grad2} is unbiased:
\[
\bbE_{Y_i\sim \pi,\, \text{i.i.d.}}[\widehat {G} ] = G(\theta)\,.
\]

\subsubsection{Coupling Method} \label{subsubsec:coupling}
Without loss of generality, we first consider the simple case where $\theta$ is a one-dimensional (1D) variable. For a fixed parameter $\theta$, we generate $N$ i.i.d.~samples, $X_1,\dots, X_N$,  following $p(x,\theta)$, and then evaluate the Monte Carlo sum $\widehat{J}(\theta)$ following~\eqref{eq:MC_obj}. Given a small perturbation $\Delta \theta$ in the parameter $\theta$, we again draw $N$  i.i.d.~samples, $\tilde X_1,\dots, \tilde  X_N$, corresponding to the perturbed density function $p(\tilde{x},\theta + \Delta \theta)$. We then calculate the new Monte Carlo sum $\widehat{J}(\theta + \Delta \theta)$. In this setup, the sampling procedure of $\{X_i\}$ and $\{\tilde X_i\}$ are independent, and we denote the resulting gradient estimation from uncorrelated distributions as $\widehat{G}_{ind}$.

Using the idea of finite-difference approximation to the first-order derivative, we can approximate the gradient $G = \nabla J(\theta) \approx \frac{1}{\Delta \theta} \left({J}(\theta + \Delta \theta) -  {J}(\theta) \right) $ using the Monte Carlo sums:
\begin{equation}\label{eq:cp_grad1}
\widehat{G}_{ind} \approx \frac{1}{\Delta \theta} \left(\widehat{J}(\theta + \Delta \theta) -  \widehat{J}(\theta) \right)=\frac{1}{N \Delta \theta} \sum_{i=1}^N \left( f(\tilde X_i) -  f(X_i)\right) \,.
\end{equation}

The variance of $\widehat{G}_{ind}$ can be calculated as follows.
\begin{eqnarray}
\text{Var}[\widehat{G}_{ind}] &=& \frac{1}{N^2 |\Delta \theta|^2} \sum_{i=1}^N 
 \left( \text{Var}_{X_i\sim p(x,\theta)} [f(X_i)] + \text{Var}_{\tilde X_i\sim p(\tilde{x},\theta+\Delta \theta)} [f(\tilde X_i)]   \right) \nonumber \\
 &=& \frac{1}{N |\Delta \theta|^2 } \left(\text{Var}_{X\sim p(x,\theta)} [f(X)] + \text{Var}_{\tilde X\sim p(\tilde{x},\theta+\Delta \theta)} [f(\tilde X)] \right)\,. \label{eq:cp_var_1}
\end{eqnarray}

Besides the random errors from computing both $\widehat{J}(\theta + \Delta \theta)$ and $\widehat{J}(\theta)$,  we address that this gradient estimator is also subject to the finite difference approximation error $\mathcal{O}(|\Delta \theta|)$. Higher-order difference schemes, such as the central difference method, can be used here. Moreover, reducing the size of $\Delta \theta$  decreases the finite difference error but, unfortunately, increases the random error as the standard deviation is proportional to $1/|\Delta \theta|$; see Equation~\eqref{eq:cp_var_1}. Therefore, to reduce the total error, $|\Delta \theta|$ cannot be chosen to be too large or too small. Certain prior analyses and estimations can be helpful to select a reasonable perturbation $\Delta \theta$; see~\cite[Sec.~6.1]{caflisch2021adjoint}.

The random error not only depends on the size of  $\Delta \theta$,  but also the concrete sampling schemes for $\{X_i\}$ and $\{\tilde X_i\}$.  The Coupling Method in the context of Monte Carlo gradients is a technique used to reduce the variance in the gradient estimation~\eqref{eq:cp_grad1} that comes from the two sets of samples. The essential idea behind the coupling method is to \textit{introduce a dependency} between pairs of random variables, 
$$
X_i, \text{ and } \tilde X_i,\quad i = 1,\dots, N\,,
$$
that are used in the Monte Carlo estimation process.

Consider a joint distribution function $P(x,\tilde{x};\theta, \Delta \theta)$ whose two marginals are $p(x,\theta)$ and $p(\tilde{x},\theta+ \Delta \theta)$,  respectively.  Let $(X_i, \tilde{X}_i)$ be i.i.d.~samples from $P(x,\tilde{x};\theta, \Delta \theta)$, $i= 1,\ldots,N$. 
We still use Equation~\eqref{eq:cp_grad1} to estimate the gradient, which is denoted by $\widehat{G}_{cp}$ to highlight the coupled sampling.   We proceed to calculate the variance of $\widehat{G}_{cp}$: 
\begin{eqnarray*}
\text{Var}_{(X_i, \tilde{X}_i) \sim P}[\widehat{G}_{cp}] &=& \frac{1}{N^2 |\Delta \theta|^2} \sum_{i=1}^N 
 \Big( \text{Var}_{X_i\sim p(x,\theta)} [f(X_i)] + \text{Var}_{\tilde X_i\sim p(\tilde,\theta+\Delta \theta)} [f(\tilde X_i)]  \\
&&\hspace{2.5cm} - 2 \text{Cov}_{(X_i, \tilde{X}_i) \sim P} [f(X_i), f(\tilde{X}_i)] \Big) \\
&=& \text{Var}[\widehat{G}_{ind}]   -  \frac{2 }{N  |\Delta \theta|^2} \text{Cov}_{(X, \tilde{X}) \sim P} [f(X), f(\tilde{X})] \,.
\end{eqnarray*}
Thus, if one can strategically correlate these two sets of random samples such that $\text{Cov}_{(X, \tilde{X}) \sim P} [f(X), f(\tilde{X})]$ is as large as possible,   the variance of the gradient estimate $\widehat{G}_{cp}$ can be significantly smaller than $\widehat{G}_{ind}$ that is obtained without coupling the two sets of samples. In other words,  we want to design correlations so that pairs of samples move together more closely in the sample space, thereby reducing the randomness in the gradient estimation in~\eqref{eq:cp_grad1}. A simple and practical way to correlate these two sets of samples is to fix the random seed in the random number generator when sampling $\{X_i\}$ and  $\{\tilde{X}_i\}$. Other coupling methods exist, and their concrete forms highly depend on the sampling algorithm under consideration.

By reducing the variance of the gradient estimates, the coupling method can lead to more efficient parameter optimization, as it requires fewer samples to achieve a certain level of accuracy. However, it is worth noting that the coupling method does not scale well when the dimension of the parameter $\theta$ increases. Consider $\theta \in \bbR^m$. Then, we need to repeat the procedure in~\eqref{eq:cp_grad1} at least $m$ times to obtain the full gradient, as each evaluation following~\eqref{eq:cp_grad1}  only yields a directional derivative. This involves $m$ different sampling tasks for one gradient evaluation. In contrast, the Score Function Method (Section~\ref{subsubsec:score}) and the Pathwise Derivative Method (Section~\ref{subsubsec:pw_grad}) require as few as one sampling process to produce a good estimation for the gradient $G$, especially if combined with the adjoint-state method to be discussed in Section~\ref{subsec:adjoint}.

We remark that the Coupling Method can be used to verify the accuracy of the gradient obtained through the adjoint Monte Carlo method~\cite{caflisch2021adjoint}.  This is because it only requires an efficient sampling algorithm, e.g., forward Monte Carlo solvers for kinetic PDEs,  for the target distribution, and does not demand access to the target density $p(x,\theta)$ or the distribution to have a reparameterization structure,  which is very hard to ensure when working with applications of kinetic PDEs.

\section{PDE-Constrained Optimization and the Adjoint-State Method}\label{sec:PDE-constrained}

Many real-world phenomena are described by PDEs, from fluid dynamics and heat transfer to electromagnetism and elasticity~\cite{borzi2011computational}. For example, as the focus of this survey, kinetic equations, which describe the evolution of systems based on their individual particle dynamics, are fundamental in various scientific and engineering disciplines. 
Optimizing systems described by these phenomena often naturally leads to PDE-constrained optimization problems. Examples include
\begin{itemize}
    \item optimizing the shape of a plasma fusion reactor,
    \item reconstructing the heat conductivity for a specific material or
    \item controlling the movements of robots.
\end{itemize}

Most PDE-constrained optimization problems share a common mathematical formulation:
\begin{equation}\label{eq:PDE-OPT}
    \min_{\theta} J(f,\theta)\quad \text{subject to}\quad h(f,\theta) = 0.
\end{equation}
Here, $J(f,\theta)$ is an objective functional of the state variable $f$ (often the PDE solution) and parameter function $\theta$. The dependence of  $f$ on $\theta$ is expressed implicitly through the operator $h$ denoting the PDE. Recall that a  PDE is an implicit formulation of its corresponding PDE solution. For convenience, we will also denote  $f = F(\theta)$ to highlight the dependence between the PDE solution $f$ and the parameter $\theta$. The (nonlinear) operator $F$ is often referred to as the forward operator in inverse problems and the control theory, and the task of solving $f$ given $h$ and $\theta$ is regarded as the forward problem.

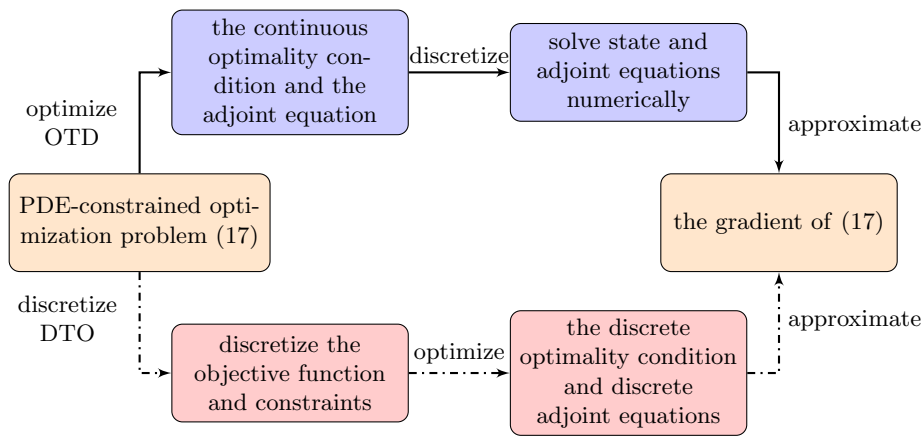
\begin{figure}
\begin{center}
{\scriptsize{
\tikzstyle{decision1} = [rectangle, draw, fill=blue!20,
    text width=6em, text centered, node distance=3cm, rounded corners, minimum height=4em]
\tikzstyle{decision2} = [rectangle, draw, fill=red!20,
    text width=6em, text centered, node distance=3cm, rounded corners, minimum height=4em]
\tikzstyle{block} = [rectangle, draw, fill=orange!20,
    text width=5em, text centered, node distance=3cm, rounded corners, minimum height=4em]
\tikzstyle{block1} = [rectangle, draw, fill=white!20,
    text width=5em, text centered, node distance=3cm, rounded corners, minimum height=4em]
\tikzstyle{line} = [draw, thick, color=black, -latex']
\tikzstyle{line2} = [draw, thick,dash dot, color=black, -latex']
\tikzstyle{line3} = [draw, thick,dash dot, color=red, -latex']
\begin{tikzpicture}[scale=1.0, node distance = 2cm, auto]
    \node(OTD1) [block,text width=10em,node distance=1cm]  {PDE-constrained optimization problem~\eqref{eq:PDE-OPT}} ;
    \node (OTD2) [decision1, above of=OTD1,text width=9em,node distance=2cm, xshift = 2cm]  {the continuous optimality condition and the adjoint equation};
    \node (DTO1) [decision2, below of=OTD1,xshift = 2cm, node distance=2cm, text width=9em]  {discretize the objective function and constraints};
    \node (OTD3) [decision1, right of=OTD2,node distance=4.5cm,text width=9em] {solve state and adjoint equations numerically};
    \node (DTO2) [decision2, right of=DTO1, node distance=4.5cm,text width=9em] {the discrete optimality condition and discrete adjoint equations};
        \node (grad) [block, above of=DTO2, xshift = 2cm, node distance=2cm,text width=9em] {the gradient of~\eqref{eq:PDE-OPT}};
\path [line] (OTD1) |- node [anchor=east,yshift = -0.7cm]{
\begin{tabular}{cc}
optimize \\
OTD\\
\end{tabular}
}(OTD2);
\path [line] (OTD2) -- node{discretize}(OTD3);
\path [line] (OTD3) -| node[anchor=west,yshift = -0.7cm]{approximate}(grad);
\path [line2] (OTD1) |-  node [anchor=east,yshift = 0.7cm]{
\begin{tabular}{cc}
discretize \\
DTO\\
\end{tabular}
}(DTO1);
\path [line2] (DTO1) -- node{optimize}(DTO2);
\path [line2] (DTO2) -| node[anchor=west,yshift = 0.7cm]{approximate}(grad);
\end{tikzpicture}
}}
\end{center}
\caption{The OTD approach (solid line) and the DTO approach (dash line) to compute the gradient with respect to the unknown parameter for a general PDE-constrained optimization problem~\eqref{eq:PDE-OPT}.}\label{fig:flowchart}
\end{figure}

\subsection{OTD and DTO Approaches}\label{subsec:OTD_DTO}
When it comes to solving~\eqref{eq:PDE-OPT},
``Optimize-then-Discretize'' (OTD) and ``Discretize-then-Optimize'' (DTO) are two fundamental approaches in the context of optimization problems constrained by PDEs. These approaches refer to the sequence in which the optimization problem is tackled with respect to its continuous formulation and discretized counterpart. For the OTD approach, the optimization problem is formulated and solved in its continuous setting. After obtaining the continuous solution (or, more commonly, its optimality condition), the governing PDEs are discretized (e.g., using finite elements, finite differences, and particle methods). For the DTO approach, the governing PDEs and the objective functional are first discretized using proper schemes, which leads to optimization problems with the finite-dimensional unknown and constraints. The optimization is then performed on this discretized problem. We outline these two approaches in a diagram; see~\Cref{fig:flowchart}.

There have been long-standing discussions on which approach is better~\cite{hinze2012discretization}. Both OTD and DTO have pros and cons, and the conclusion highly depends on the specific PDE-constrained optimization problem under consideration and the choice of numerical scheme for discretization. OTD often leads to well-structured optimality conditions and can provide insights into the nature of the optimal solution in terms of analysis. The resulting discrete problem might inherit only some of the properties of the continuous problem. One has to conduct discretization with caution to avoid inconsistency issues that can lead to lower-order methods or even a wrong solution. For the DTO approach,  the optimization is performed on a finite-dimensional space, which can be more straightforward and amenable to standard optimization techniques. This approach is often more practical for numerical implementations. The discrete problem might preserve only some of the properties of the continuous problem, as we will see later in Sections~\ref{sec:rte}-\ref{sec:Boltzmann}. The resulting discrete optimization problem can also be large-scale, especially for fine discretizations. The DTO approach is particularly preferred for problems that are otherwise unsolvable under the OTD approach, such as the heat-transfer optimization problem~\cite{betts2005discretize,ghobadi2009discretize}. The OTD approach is very attractive if one wants to utilize the rich geometric structure of the continuous function space~\cite{jacobs2019solving}.

For a given discretization scheme, the gradient obtained from  OTD and DTO may not be equivalent, as demonstrated by many examples in the literature~\cite{hager2000runge,burkardt2002insensitive,
abraham2004effect,hinze2012discretization,liu2019non}. The lack of symmetry in discretization often causes such inconsistency. Examples include the non-Galerkin discretization of the optimality system. Such differences may contribute to (big) errors in optimization for solving nonlinear problems. 

\subsection{Adjoint-State Method}\label{subsec:adjoint}
Both OTD and DTO approaches require deriving the first-order optimality condition, especially the derivative of the (approximated) objective $J$ with respect to the (discretized) parameter $\theta$. Next, we introduce an efficient way to compute the  derivative. Under the OTD approach, $f$ and  $\theta$ are functions while in the DTO approach, they are both discrete objects as a result of discretization. In the following derivations, we will not distinguish between these two.

Since the PDE constraint $h$ implies that $f$ depends on $\theta$, we will plug in $f = F(\theta)$ into~\eqref{eq:PDE-OPT}. Thus, we deal with the optimization problem
\begin{equation}\label{eq:PDE-OPT-2}
    \min_{\theta} \mathfrak{J}(\theta)\,, \text{\quad where \quad }\mathfrak{J}(\theta)  = J(F(\theta) , \theta)\,.
\end{equation}
Correspondingly, the PDE constraint becomes
\begin{equation}\label{eq:constraint}
    h(F(\theta),\theta) = 0.
\end{equation}
 Without loss of generality, we assign all the variables with corresponding linear spaces that can be both finite and infinite dimensional:  $f \in U$,  $\theta \in V$,  $h(f,u) \in W$ and $U,V,W$ are Hilbert spaces. We further assume that operators $h$ and $F$ are smooth enough such that they are Fr\'echet differentiable.

A first-order variation of $h$ with respect to $\theta$ yields
\begin{equation} \label{eq:linear_h}
   h_f (F(\theta), \theta) \, D_\theta F(\theta) +   h_\theta (F(\theta), \theta)  = 0 \, ,
\end{equation}
where $h_f$ and $h_\theta$ are Fr\'echet derivatives of $h(f,\theta)$ with respect to $f$ and $\theta$ variables, and $D_\theta F$ is the derivative of $F(\theta)$ with respect to $\theta$. Note that $h_f \in \mathcal{L}(U,W)$, $h_\theta \in \mathcal{L}(V,W)$ and $D_\theta F \in \mathcal{L}(V,U)$, where $\mathcal{L}(A,B)$ denotes the linear operator that maps elements  from function space $A$ to function space $B$.

To calculate the derivative of $\mathfrak{J}$ with respect to $\theta$, we use chain rule and obtain
\begin{equation}\label{eq:grad_1}
 D_\theta \mathfrak{J} =   J_f(f, \theta)  D_\theta F(\theta) + J_\theta(f,\theta)\,,
\end{equation}
where $D_\theta\mathfrak{J} \in V^*$, the dual space of $V$. Also,  with $U^*$ being the the dual space of $U$, $ J_f \in U^*$  and $J_\theta \in V^* $ are partial Fr\'echet derivatives of $J(f,\theta)$ with respect to $f$ and $\theta$ variables, respectively.  

Using~\eqref{eq:grad_1}  to compute the derivative of $\mathfrak{J}$ requires  $D_\theta F$, the linearized operator of the forward map $F$. Recall that the dependence between $\theta$ and $f$ is \textit{implicitly} given through the PDE constraint $h$. That is, $F$ is defined through~\eqref{eq:constraint}. As a result, to obtain $D_\theta F$,  we have to utilize~\eqref{eq:linear_h}. One may first solve for $D_\theta F(\theta)$ at a particular $\theta^*$ 
through
\begin{equation}\label{eq:jac_solve}
h_f (F(\theta^*), \theta^*) \,  D_\theta F (\theta^*) = - h_\theta(F(\theta^*), \theta^*)\,,
\end{equation}
and then plug the solution into~\eqref{eq:grad_1}. We remark that since the objective function $J$ is user's choice, $J_f$ and $J_\theta$ are usually known and easy to compute based on the form of $J$. If $V$, the space of $\theta$, is finite-dimensional, then the cost of this inverse solve~\eqref{eq:jac_solve} increases linearly with respect to the dimension of $V$ (i.e., the intrinsic degrees of freedom in the unknown parameter $\theta$), and can be computationally prohibitive for very large-scale problems.

Next, we present a more efficient way of obtaining $D_\theta \mathfrak{J}$, called the \textit{Adjoint-State Method}~\cite{hinze2008optimization,borzi2011computational} or the reduced gradient method~\cite{herzog2010algorithms} in the literature. We first derive it using the method of Lagrangian multipliers. Consider the Lagrangian $L$ given by
\begin{equation}\label{eq:PDE-LAG}
    L(f,\lambda,\theta) = J(f,\theta) +  \lambda  h(f,\theta) ,
\end{equation}
where $\lambda \in W^*$, the dual space of $W$. In the Lagrangian~\eqref{eq:PDE-LAG}, there is no dependence among $f,\lambda$ and $\theta$, as all of them are independent variables. 
However, if we plug $f = F(\theta)$ into~\eqref{eq:PDE-LAG}, which is equivalent to restraining the Lagrangian over a particular manifold $\mathcal{M} := \{(f,\theta): f = F(\theta)\}$, a subset of the full space $U \times V$, 
we would obtain
\begin{eqnarray}
\mathfrak{L}(\lambda,\theta)&:=&L(F(\theta),\lambda,\theta) \nonumber \\
&=& J(F(\theta),\theta) + \lambda  h(F(\theta),\theta) \nonumber \\
&=&  J(F(\theta),\theta) \nonumber \\
&=& \mathfrak{J}(\theta), \label{eq:PDE-LAG-2}
\end{eqnarray}
where we used~\eqref{eq:PDE-OPT-2} and~\eqref{eq:constraint}. Furthermore, the derivative of $\mathfrak{L}$ with respect to $\theta$ yields
\begin{eqnarray}\label{eq:PDE-LAG-3}
     \mathfrak{L}_\theta( \lambda,\theta)   = L_\theta (F(\theta), \lambda,\theta) + L_f(F(\theta), \lambda,\theta) \, D_\theta F(\theta)\,,
\end{eqnarray}
where $\mathfrak{L}_\theta, L_\theta \in V^*$ and $L_f \in U^*$. We remark that both $\mathfrak{L}$ and $ \mathfrak{L}_\theta$ are independent of the variable $\lambda$ since the last term in~\eqref{eq:PDE-LAG-2} is $\lambda$-independent.

If we choose a particular $\lambda \in W^*$ such that $L_f(F(\theta), \lambda,\theta) = 0 \in U^*$, then we do not have to involve $D_\theta F(\theta)$ for the calculation of the derivative following the formula~\eqref{eq:PDE-LAG-3}.  \textbf{This is the central idea behind the adjoint-state method.} After further simplification for $L_f(F(\theta), \lambda,\theta)$, we obtain
\[
    L_f(F(\theta), \lambda,\theta) = J_f (F(\theta) ,\theta) +   \lambda \, h_f(F(\theta) ,\theta))\,.
\]
If $\lambda$ solves the following equation 
\begin{equation}\label{eq:adj_eqn}
   \lambda \, h_f(F(\theta) ,\theta)) = - J_f (F(\theta) ,\theta),
\end{equation}
then $\lambda$ is $\theta$-dependent and we denote it by $\lambda(\theta)$. Finally, with this particular choice of $\lambda$, we have reached a simple form of the derivative following~\eqref{eq:PDE-LAG-3}:
\begin{eqnarray}
    D_\theta \mathfrak{J}(\theta) & = & D_\theta \mathfrak{L}(F(\theta), \lambda(\theta),\theta) \nonumber \\
    &=& L_\theta (F(\theta), \lambda(\theta),\theta) \nonumber \\
    &=& J_\theta(F(\theta) ,\theta) + \lambda(\theta)\, h_\theta(F(\theta) ,\theta)). \label{eq:derivative}
\end{eqnarray}
The first equation holds because of~\eqref{eq:PDE-LAG-2}.

Equation~\eqref{eq:adj_eqn} is often referred to as the \textit{adjoint equation}, while solving $f$ based on~\eqref{eq:constraint} is referred to as the \textit{state equation}. In the \textit{adjoint-state method}, for a given $\theta^*$, we solve the forward equation once and obtain $f = F(\theta^*)$. It is followed by an adjoint equation solve, which results in $\lambda$. Note that $\lambda$ implicitly  depends on $\theta^*$ due to Equation~\eqref{eq:adj_eqn}. We then obtain the derivative $D_\theta \mathfrak{J}(\theta^*)$ of the constrained problem following the formula~\eqref{eq:derivative}. If $h$ is a PDE constraint, this amounts to two PDE solves, independent of the dimension of the discretized $\theta$-space $V$, making it extremely attractive for infinite-dimensional optimization problems. Therefore, the adjoint-state method is widely used for large-scale practical problems with a differentiable equality constraint~\eqref{eq:constraint} and a differentiable objective functional $J$.

\begin{rem}
We have presented a framework that can apply to both OTD and DTO approaches for a large class of PDE-constrained problems. Note that this requires $U$ (the space for $f$), $V$ (the space for $\theta$), and $W$ (the space for the range of the equality constraint $h$) to be vector spaces and, preferably, Hilbert spaces. However, for many kinetic equations, the solution should be regarded as (probability) distribution functions with nonnegativity and total mass conserved. Moreover, for discretization schemes based on particle methods, the relationship between the continuous and discrete variables is not linear, not to mention the further complications originating from the stochasticity in random algorithms (e.g., the Monte Carlo method). In the later sections, we will discuss two specific cases and how to \underline{adapt} the adjoint-state method to such representative kinetic equation-based optimization problems. 
\end{rem}
 
In the following sections, we will study closely two kinetic equation-constrained optimization problems and how to employ and adapt the general framework in Section~\ref{sec:PDE-constrained} to solve them systemically. Since both the radiative transport equation and the Boltzmann equation can be solved through Monte Carlo particle methods, we will also derive the corresponding adjoint Monte Carlo method for each of the two kinetic equations. Interested readers may refer to~\cite{caflisch2021adjoint,li2022monte,yang2023adjoint} for more details in terms of motivation, derivation and convergence proofs.

\section{Case Study I: Radiative Transport Equation}\label{sec:rte}
The radiative transport equation (RTE), also known as the radiative transfer equation, is a model problem for simulating light propagation in an optical environment. For  simplicity, we consider time-dependent RTE with no spatial boundary effect,
\begin{equation}\label{eqn:RTE}
\begin{cases}
\partial_t f +v\cdot\nabla_xf &= \, \sigma(x)\mathcal{L}[f], \\
\quad f(t=0,x,v) &= \, f_0 (x,v),
\end{cases}\qquad x\in \mathbb{R}^{d_x}, ~ v \in \Omega\,.
\end{equation}
Here, $f(t,x,v)$ is the distribution function of photon particles at time $t$ on the phase space $(x,v)$ and $d_x$ is the dimension for the spatial domain. The left-hand side of the equation describes the photon moving in a straight line in $x$ with velocity $v$, whereas the right-hand side characterizes the photon particles' interaction with the media characterized by the function $\sigma(x)$. The term $\mathcal{L}[f]$ is defined by
\begin{equation}\label{eq:f_v}
\mathcal{L}[f] = \frac{1}{|\Omega|} \langle f\rangle_v - f\,,\quad\text{with}\quad \langle f\rangle_v =  \int_\Omega f(x,v)\rd{v}\,.
\end{equation}
The entire right-hand side represents that particles at location $x$ have a probability proportional to $\sigma(x)$ to be scattered, into a new direction uniformly chosen in the $v$ space. As a result, the distribution function in the phase space exhibits a gain on $|\Omega|^{-1}\langle f\rangle_v(x)$ and a loss on $f(x,v)$. We will use $\langle\cdot\rangle_s$ to denote the integration in variable $s$ with respect to the Lebesgue measure.

The forward problem~\eqref{eqn:RTE}, i.e., given $\sigma$ and solving for $f$, has a unique solution under very mild conditions on both $\sigma(x)$ and the initial data $f_0$~\cite{dautray1993mathematical,bal2006radiative}.
Moreover, the equation preserves mass because $\rho=\langle f\rangle_{xv}$ is a constant in time. Without loss of generality, we set $\rho=1$. Note that~\eqref{eqn:RTE} is a linear equation with respect to $f$.

\subsection{Constrained Optimization Problem}
We consider a PDE-constrained optimization originated from an inverse problem as an example. Assume we have measurement which is the final-time solution of the RTE,
\[
d(x,v) \approx f(t=T,x,v)\,,
\] 
and the goal is to identify the variable absorption parameter $\sigma(x)$ in the spatial domain given the known initial condition. 
The objective functional for the corresponding computational inverse problem is the mismatch of the simulated solution and the measured data, and the constraints come from the fact that the simulated solution  satisfies a forward RTE. The corresponding PDE-constrained optimization is of the form
\begin{equation}\label{eqn:RTE_min}
\min_\sigma J\left(f(\boldsymbol{\cdot}\, ; \sigma)\right), \quad \text{s.t.} \text{~$f(\boldsymbol{\cdot}\, ; \sigma)$ satisfies~\eqref{eqn:RTE}\,,} 
\end{equation}
where $f(\boldsymbol{\cdot}\,;\sigma)$ is the simulated data that solves~\eqref{eqn:RTE} with the given initial condition and the absorbing parameter $\sigma(x)$.
In this particular application setup, we may consider the objective functional
\[
J\left(f(\boldsymbol{\cdot}\, ; \sigma)\right)=\frac{1}{2}\iint  
 | f(t=T,x,v;\sigma) - d(x,v) |^2\rd x \rd v\,.
\]
Similar to~\eqref{eq:constraint}, we can rewrite~\eqref{eqn:RTE} to match the general form of an equality constraint $h(f,\sigma) = 0$. We then denote $f = F(\sigma)$ with $F$ being the forward operator. Next, we consider all relevant function spaces, $U$ (the space of $f$), $V$ (the space of $\sigma$), $W$ (the space of the constraint) to be $L^2$ over the respective domain.

\subsection{The OTD Approach}
In the following, we detail how to employ the OTD approach to  design an adjoint Monte Carlo method to tackle~\eqref{eqn:RTE_min}.
\subsubsection{Gradient Derivation}
We first apply the method of Lagrange multipliers and consider the Lagrangian following the recipe given in Section~\ref{subsec:adjoint}.
\[
L(f,\lambda, \bar{\lambda}, \sigma) = J(f) + \langle \lambda \,,\partial_t f +v\cdot\nabla_x f -  \sigma(x)\mathcal{L}[f]\rangle_{txv} + \langle \bar{\lambda}\,,f(t=0,x,v) - f_0(x,v)\rangle_{xv}.
\]
Here, functions $\lambda(t,x,v)$ and $\bar{\lambda}(x,v)$ are Lagrange multipliers with respect to the RTE solution $f(t,x,v)$ for $t>0$ and the initial condition $f(t=0,x,v)$, respectively. 

We remark that the $\lambda$ and $\bar{\lambda}$ we used here slightly differ from the adjoint variable $\lambda$ in Section~\ref{subsec:adjoint}. The former is in the function space $W$ where $h(f,\sigma)$ belongs to,  while the latter is considered an element of $W^*$, i.e., the \textit{dual space} of $W$. In other words, the $\lambda,\bar{\lambda}$ used here represent the Riesz representation of the $\lambda$ in~\eqref{eq:PDE-LAG} through the Riesz map. Since all the spaces we consider in this subsection are $L^2$ for the RTE-based constrained optimization, these two coincide, and thus, we will not differentiate between them when referring to the adjoint variable.

With integration by parts, we rewrite the Lagrangian as:
\[
\begin{aligned}
L(f,\lambda, \bar{\lambda}, \sigma) = & J(f) +  \langle f\,, -\partial_t \lambda - v\cdot\nabla_x \lambda - \sigma(x)\mathcal{L}[\lambda]\rangle_{txv} + \langle \bar{\lambda}\,,f(t=0,x,v) - f_0(x,v)\rangle_{xv}\\
& + \langle \lambda(t=T,x,v)\,, f(t=T,x,v)\rangle_{xv} - \langle \lambda(t=0,x,v)\,, f(t=0,x,v)\rangle_{xv} \,.
\end{aligned}
\]
Based on~\eqref{eq:adj_eqn}, we can set the derivative of $L(f,\lambda, \bar{\lambda}, \sigma)$ with respect to $f(t,x,v)$ for any $0\leq t \leq T$ evaluated at $f = F(\sigma)$ to be zero. We then obtain the adjoint equation for the RTE-constrained optimization problem:
\begin{equation}\label{eqn:g_eqn}
\begin{cases}
- \partial_t \lambda  - v\cdot\nabla_x \lambda &=  \sigma(x)\mathcal{L}[\lambda]\,, \\
\qquad \lambda(T,x,v) &=   - \frac{\delta J}{\delta f(T,x,v)} \,,
\end{cases}
\end{equation}
where the final condition at $T$ comes from the form of $J$ given in~\eqref{eqn:RTE_min}. Note that the resulting adjoint solution $\lambda$ is $\sigma$-dependent.

Following the general framework in~\eqref{eq:derivative}, we then have the Fr\'echet derivative of $\mathfrak{J}(\sigma) := J(F(\sigma))$ as
\begin{equation}\label{eq:RTE_OTD_grad}
    D_\sigma \mathfrak{J}(\sigma) = L_\sigma(F(\sigma),\lambda,\bar{\lambda}, \sigma) =  - \langle \lambda , \mathcal{L}[f] \rangle_{tv} = \int_0^T  \langle \lambda f\rangle_v \rd t - \frac{1}{|\Omega| } \int_{0}^T \langle f\rangle_v \langle \lambda \rangle_v \rd t.
\end{equation}

The derivation above is carried out completely on the function space in the continuous setting, as done in the OTD framework. Here, we consider the function space for the parameter $\sigma(x)$ to be $L^2$, so the derivative $D_\sigma \mathfrak{J}(\sigma)$ is indeed the $L^2$ gradient of $\mathfrak{J}$. 

\subsubsection{Discretization}
Upon obtaining \eqref{eq:RTE_OTD_grad}, the next step is discretization, which involves numerically solving (i) the state equation~\eqref{eqn:RTE}, (ii) the adjoint equation~\eqref{eqn:g_eqn}, and (iii) a discrete approximation for the gradient~\eqref{eq:RTE_OTD_grad}. 

Assume we have a good numerical scheme for the forward state equation~\eqref{eqn:RTE}. It is also quite tempting to solve the adjoint equation~\eqref{eqn:g_eqn} using the same scheme, given the great similarity between~\eqref{eqn:RTE} and~\eqref{eqn:g_eqn}, which results from the fact that the state equation~\eqref{eqn:RTE} is linear with respect to $f$. However, one has to be careful here as we want to construct a consistent numerical scheme for the adjoint equation with respect to the one applied to the forward equation. That is, the \textit{numerical}  solution for the adjoint equation~\eqref{eqn:g_eqn} is preferably to be an adjoint of the \textit{numerical} solution for the state equation~\eqref{eqn:RTE}. Failure to preserve the consistency, i.e., the adjoint relationship on the discrete level, may result in a low-order or even wrong gradient approximation~\cite{hager2000runge}.

First, we consider the finite-volume method for the state equation~\eqref{eqn:RTE}. We only consider the 1D case for both $x$ and $v$ as an illustrating example. Denote
\[
\fij^m \approx f(t_m, x_i, v_j)\,,\quad  \gij^m \approx \lambda(t_m, x_i, v_j)\,, \quad 0\leq m \leq M, \quad 1\leq j \leq N_v,\quad 1\leq i \leq N_x\,,
\]
as the numerical approximation, and let  $\dx$, $\dt$ and  $\dv $ be the corresponding mesh size in $x$, $t$ and $v$, respectively. We consider the time domain $[0,T]$, the spatial domain $D = [-2,2]$ and the velocity domain $\Omega = [-1,1]$. Then $M \dt = T$,  $t_m = m \dt$,  $\dx N_x = |D| = 4$,  and $\dv N_v = |\Omega| = 2$. We have the following discretization for \eqref{eqn:RTE}:
\begin{align}\label{eq:fvm_f}
    \fij^\np - \fij^m + \frac{\dt}{\dx} v_j^{+}(\fij^m - \fijm^m) 
    + \frac{\dt}{\dx} v_j^{-}(\fijp^m - \fij^m) = \sigma_i \dt \left(|\Omega|^{-1}\average{\fij^m}_v - \fij^m \right)
    \,,
\end{align}
for $0\leq m \leq M -1$, $1\leq i \leq N_x$,  and $1\leq j \leq N_v$, 
with the initial condition 
\begin{align*}
    \fij^0 = f_0(x_i, v_j)\,.
\end{align*}
Here, $v_j = -1+ (j-1/2) \dv$, $j = 0,\ldots, N_v$, representing the cell center, and $v^+ = \max\{v,0\}$, $v^- = \min\{v,0\}$. The average in $v$ is computed by a simple midpoint rule: $\average{\fij^n}_v = \sum_{j=1}^{N_v} \fij^n \dv$.  The consistent scheme for $\lambda$ in the adjoint equation~\eqref{eqn:g_eqn} is
\begin{align}\label{eq:fvm_lambda}
     - (\gij^\np - \gij^m ) - \frac{\dt}{\dx} v_j^{+}(\gijp^m - \gij^m ) 
    - \frac{\dt}{\dx} v_j^{-}( \gij^m - \gijm^m ) = \sigma_i \dt (|\Omega|^{-1}\average{\gij^\np}_v - \gij^\np)
    \,,
\end{align}
for $m = M-1, M-2, \ldots, \geq 0$ and $1\leq j \leq N_v$, with the final condition depending on the objective function $J$.

Note that $v_j^+$ is multiplied by $(\fij^m - \fijm^m)$ in~\eqref{eq:fvm_f}, but by $(\gijp^m - \gij^m )$ in~\eqref{eq:fvm_lambda}. The right-hand sides of~\eqref{eq:fvm_f} and~\eqref{eq:fvm_lambda} are also evaluated at different time points.  In other words, we cannot simply replace $\{\fij^m\}$  in~\eqref{eq:fvm_f}
 with $\{\gij^m\}$ and flip the sign on the right-hand side as a numerical solver for the adjoint equation~\eqref{eqn:g_eqn}, despite the fact that the continuous adjoint equation~\eqref{eqn:g_eqn} can be derived by applying these steps to the continuous forward equation~\eqref{eqn:RTE}.  These are the differences between consistent forward and adjoint numerical schemes and inconsistent ones.
 
\subsubsection{Discretization with Adjoint Particle Method}\label{subsec:P-OTD}
The RTE in~\eqref{eqn:RTE} has an efficient Monte Carlo solver with provable convergence~\cite[Theorem 2.1]{li2022monte}. The step-by-step numerical scheme is detailed in Algorithm~\ref{alg:f-RTE}. 
\begin{algorithm}
\caption{Monte Carlo Method for Solving the Forward RTE~\eqref{eqn:RTE}\label{alg:f-RTE}}
\begin{algorithmic}[1]
\State  \textbf{Preparation:} $N$ pairs of samples $\{(x_n^0,v_n^0)\}_{n=1}^N$ from the initial distribution $f_0(x,v)$; the total time steps $M$ and the time spacing $\Delta t$; and the parameter function $\sigma(x)$.
\For{$m=0$ to $M-1$}
\State Given $\{(x_n^m,v_n^m)\}_{n=1}^N$, set $x_n^{m+1} = x_n^m + \Delta t \,  v_n^m$, $n=1,\ldots,N$.
\State Draw random numbers $\{p^{m+1}_n\}_{n=1}^N$ from the uniform distribution $\mathcal{U}([0,1])$.
\If{ $p_n^{m+1} \geq \alpha_n^{m+1} = \exp(-\sigma(x_n^{m+1}) \Delta t)$}  $v_n^{m+1} = \eta_n^{m+1}$ where $\eta_n^{m+1}\sim\mathcal{U}(\Omega)$.
\Else  \, $v_n^{m+1} = v_n^{m}$.
\EndIf
\EndFor
\end{algorithmic}
\end{algorithm}
As noted earlier, the adjoint equation~\eqref{eqn:g_eqn} can be seen as the forward equation~\eqref{eqn:RTE} with the sign in time and velocity reversed, which implies that we can also adapt the same Monte Carlo method for solving the adjoint equation.  If we follow this approach, an immediate  problem will arise as detailed in~\cite[Section 2.2]{li2022monte}: both $f$ and $\lambda$ are approximated by empirical distributions, i.e., finite sums of delta measures, and consequently the gradient evaluation in~\eqref{eq:RTE_OTD_grad} will involve the multiplication of two delta measures over the $x$ and $v$ domains, and these two sets of trajectories meet in the $(t,x,v)$ space with probability zero. 
Therefore, a direct discretization of the adjoint equation~\eqref{eqn:g_eqn} using the same Monte Carlo method leads to \textit{severe inconsistency} between the forward and adjoint numerical schemes.

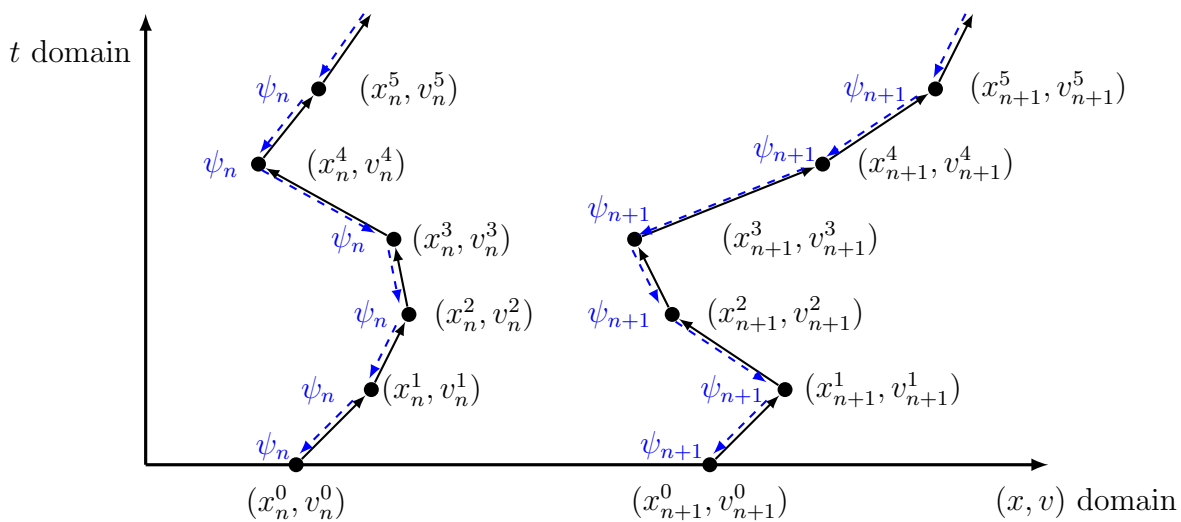
\begin{figure}
\begin{center}
    \begin{tikzpicture}
    \draw [-latex, very thick] (-6,0) -- (-6,6);
    \draw [-latex, very thick] (-6,0) -- (6,0);
    \node at (6.5,-0.5) {$(x,v)$ domain};
    \node at (-7,5.5) {$t$ domain};
    \node(x0) [circle,fill,inner sep=2pt] at (-4,0) {};
    \node(x1) [circle,fill,inner sep=2pt] at (-3,1) {};
    \node(x2) [circle,fill,inner sep=2pt] at (-2.5,2) {};
    \node(x3) [circle,fill,inner sep=2pt] at (-2.7,3) {};
    \node(x4) [circle,fill,inner sep=2pt] at (-4.5,4) {};
    \node(x5) [circle,fill,inner sep=2pt] at (-3.7,5) {};
    \node(xx0) at (-4.1,0) {};
    \node(xx1) at (-3.1,1) {};
    \node(xx2) at (-2.6,2) {};
    \node(xx3) at (-2.8,3) {};
    \node(xx4)  at (-4.6,4) {};
    \node(xx5) at (-3.8,5) {};
    \node at (-4, -0.5) {$(x_n^0, v_n^0)$} ;
     \node at (-2.2, 1) {$(x_n^1, v_n^1)$} ;
    \node at (-1.5, 2) {$(x_n^2, v_n^2)$} ;
    \node at (-1.8, 3) {$(x_n^3, v_n^3)$} ;
    \node at (-3.2, 4) {$(x_n^4, v_n^4)$} ;
    \node at (-2.5, 5) {$(x_n^5, v_n^5)$} ;
    \node at (-4.3, 0.25) [blue] {$\psi_n$} ;
     \node at (-3.7, 1) [blue]  {$\psi_n$} ;
    \node at (-3, 2) [blue]  {$\psi_n$} ;
    \node at (-3.3, 3) [blue]  {$\psi_n$};
    \node at (-5, 4) [blue]  {$\psi_n$} ;
    \node at (-4.3, 5) [blue]  {$\psi_n$} ;
    \draw [thick,-latex] (x0) -- (x1);
    \draw [thick,-latex] (x1) -- (x2);
    \draw [thick,-latex] (x2) -- (x3);
    \draw [thick,-latex] (x3) -- (x4);
    \draw [thick,-latex] (x4) -- (x5);
    \draw [thick,-latex] (x5) -- (-3,6);
    \draw [blue,dashed,thick,-latex] (xx1) -- (xx0);
    \draw [blue, dashed,thick,-latex] (xx2) -- (xx1);
    \draw [blue,dashed, thick,-latex] (xx3) -- (xx2);
    \draw [blue,dashed,thick,-latex] (xx4) -- (xx3);
    \draw [blue,dashed, thick,-latex] (xx5) -- (xx4);
    \draw [blue, dashed, thick,-latex] (-3.1,6) -- (xx5);
    \node(y0) [circle,fill,inner sep=2pt] at (1.5,0) {};
    \node(y1) [circle,fill,inner sep=2pt] at (2.5,1) {};
    \node(y2) [circle,fill,inner sep=2pt] at (1,2) {};
    \node(y3) [circle,fill,inner sep=2pt] at (0.5,3) {};
    \node(y4) [circle,fill,inner sep=2pt] at (3,4) {};
    \node(y5) [circle,fill,inner sep=2pt] at (4.5,5) {};
    \node(yy0) at (1.4,0) {};
    \node(yy1) at (2.4,1) {};
    \node(yy2) at (0.9,2) {};
    \node(yy3) at (0.4,3) {};
    \node(yy4) at (2.9,4) {};
    \node(yy5) at (4.4,5) {};
    \node at (2.5, 2) {$(x_{n+1}^2, v_{n+1}^2)$} ;
    \node at (2.7, 3) {$(x_{n+1}^3, v_{n+1}^3)$} ;
    \node at (4.5, 4) {$(x_{n+1}^4, v_{n+1}^4)$} ;
    \node at (6, 5) {$(x_{n+1}^5, v_{n+1}^5)$} ;
    \node at (3.8, 1) {$(x_{n+1}^1, v_{n+1}^1)$} ;
    \node at (1.5, -0.5){$(x_{n+1}^0, v_{n+1}^0)$} ;
    \node at (0.3, 2) [blue]{$\psi_{n+1}$} ;
    \node at (0.3, 3.4) [blue]{$\psi_{n+1}$} ;
    \node at (2.5, 4.2) [blue]{$\psi_{n+1}$} ;
    \node at (3.7, 5) [blue]{$\psi_{n+1}$} ;
    \node at (1.8, 1) [blue]{$\psi_{n+1}$} ;
    \node at (1, 0.25)[blue]{$\psi_{n+1}$} ;
    \draw [thick,-latex] (y0) -- (y1);
    \draw [thick,-latex] (y1) -- (y2);
    \draw [thick,-latex] (y2) -- (y3);
    \draw [thick,-latex] (y3) -- (y4);
    \draw [thick,-latex] (y4) -- (y5);
    \draw [thick,-latex] (y5) -- (5,6);
    \draw [blue,dashed,thick,-latex] (yy1) -- (yy0);
    \draw [blue,dashed,thick,-latex] (yy2) -- (yy1);
    \draw [blue,dashed,thick,-latex] (yy3) -- (yy2);
    \draw [blue,dashed,thick,-latex] (yy4) -- (yy3);
    \draw [blue,dashed,thick,-latex] (yy5) -- (yy4);
    \draw [blue,dashed,thick,-latex] (4.9,6) -- (yy5);
    \end{tikzpicture}
\end{center}
\caption{An illustration of the RTE forward Monte Carlo trajectory (black) and its adjoint Monte Carlo trajectory (blue) for particle indices $n$ and $n+1$, as examples. The forward and adjoint trajectories coincide since they share the same randomness in the adjoint Monte Carlo method. Moreover, the adjoint variable $\gamma(t,x,v)$ carries a fixed value $\psi_n$ (resp.~$\psi_{n+1}$) along the time trajectory; see~\eqref{eq:rte_adj_MC}.\label{fig:RTE}}
\end{figure}

To circumvent this issue, \cite{li2022monte} proposed a different adjoint Monte Carlo approach, which is also much more simple than~\Cref{alg:f-RTE}. First, we observe that, on the continuous level, 
\[
\partial_t \average{fg}_{xv} \equiv 0.
\]
Since $f$ is approximated by an empirical distribution 
\begin{equation}\label{eq:f_N}
f_N(t,x,v) \approx \frac{1}{N}\sum_{n=1}^N \delta\left(x- {x_n(t)} \right) \delta\left(v- {v_n(t)} \right)\,,
\end{equation}
by combining the previous equations we obtain the relation that
\[
    \sum_{n=1}^N \lambda(t_m, x_n^m, v_n^m) =    \sum_{n=1}^N \lambda( t_{m+1},x_n^{m+1}, v_n^{m+1})\,.
\]
It is sufficient to require that for $m= M-1,\ldots,0$
\[
\lambda(t_m, x_n^m, v_n^m)  = \lambda(t_{m+1}, x_n^{m+1}, v_n^{m+1}) =\cdots =\lambda(t^{M}, x_n^{M}, v_n^{M}) = -\frac{\delta J}{\delta f}(T,x_n^M,v_n^M)=:\psi_n \,.
\]
The final-time condition comes from~\eqref{eqn:g_eqn}. In other words, we do not even need to ``solve'' the adjoint equation~\eqref{eqn:g_eqn}, as $\lambda(t,x,v)$  at time $t=t_m$ is represented by $N$ weighted adjoint particles
\begin{equation}\label{eq:rte_adj_MC}
\{\left( x=x_n^m,\, v=v_n^m, \, \psi_n \right)\}_{n=1}^N\,.
\end{equation}
Each particle is located at $x=x_n^m$ and  $v=v_n^m$, with weight $\psi_n$.
We can perform numerical interpolation and obtain a function approximating $\lambda(t,x,v)$ with certain assumptions on the smoothness of $\lambda(t,x,v)$. We plot two trajectories with indices $n$ and $n+1$ in~\Cref{fig:RTE} as examples, illustrating that the adjoint variable carries the same value $\psi_n$ and $\psi_{n+1}$ propagating backward in time along the same trajectory from the forward Monte Carlo solver.

Moreover, we have obtained all values  in $\{ \left( x_n^m,\, v_n^m, \, \psi_n \right) \}$ immediately after solving the forward equation~\eqref{eqn:RTE} using the Monte Carlo Algorithm~\ref{alg:f-RTE}. In contrast to the finite volume scheme mentioned above, this adjoint Monte Carlo approach requires zero computational cost in backward time evolution but could potentially demand flops if a continuous format of $\lambda(t,x,v)$ is needed.

The approximation for the gradient~\eqref{eq:RTE_OTD_grad} by using the Monte Carlo solver for $f$ and the weighted particle representation~\eqref{eq:rte_adj_MC} for $\lambda$ is proved to converge to the true gradient; see~\cite[Theorem 2.2]{li2022monte} for details. In this step, we fully \textit{correlate} the forward and adjoint Monte Carlo solves in terms of their source of randomness. This is analogous to the \textbf{Coupling Method} discussed in Section~\ref{subsubsec:coupling}. Unlike the Monte Carlo gradient estimator from the standard coupling method, the gradient obtained here  does not suffer from the finite difference error, and more importantly, the computational cost is \textbf{independent} of the size of the discretized $\theta$.

This is our first example of using the adjoint Monte Carlo method as an efficient approach for kinetic equation-constrained optimization problems.

\subsection{The DTO Approach}\label{subsec:RTE_DTO}
In this subsection, we will employ an alternative approach, the DTO framework, to compute the gradient of the RTE-constrained optimization problem~\eqref{eqn:RTE_min} through a \textit{different} adjoint Monte Carlo method. As the first step for the DTO approach, we need to discretize the objective functional and the state equation~\eqref{eqn:RTE}. If we use the finite-volume method, the discrete state equation would be~\eqref{eq:fvm_f}, and the corresponding discrete adjoint equation follows~\eqref{eq:fvm_lambda}. Interested reader may refer to~\cite[Appendix A]{li2022monte} for details. In the following, we consider the discretization based on the Monte Carlo Algorithm~\ref{alg:f-RTE}.

Without loss of generality, we consider the objective functional
\begin{equation}\label{eq:J2}
    J = \iint r(x,v) f(T,x,v) \rd x \rd v  = \mathbb{E}_{(X,V)\sim  f(T,x,v) } \left[r(X,V)\right]\,,
\end{equation}
where $J$ denotes the expectation of $r$ with respect to the  final-time RTE solution $f(T,x,v)$. The choice of $J$ is only for notational convenience, and the derivations shall easily apply to general objective functionals. If $f$ is approximated by $f_N$ given in Equation~\eqref{eq:f_N}, the value in \eqref{eq:J2} could be approximated by the Monte Carlo quadrature
\begin{equation}\label{eq:J2_MC}
J \approx \frac{1}{N} \sum_{n=1}^N r\left(x_n^M, v_n^M \right) =: \hat{J}  \,.
\end{equation}
Next, we assume that for any $m= 0,\ldots M$, 
$$
(x_n^m, v_n^m) \,\, \overset{\text{i.i.d.}} {\sim} \,\, f(t_m,x,v),\quad \forall n=1,\ldots,N\,.
$$
 
Given any smooth test function $\phi(v)$ and $m = 0,\ldots, M-1$, we have
\begin{equation} \label{eq:probab}
    \bbE_{\eta^{m+1}} \bbE_{p^{m+1}} [\phi(v) | (x_n^m,v_n^m) ] = \alpha_n^{m+1}  \phi(v_n^m) + (1-\alpha_n^{m+1} ) \frac{1}{|\Omega|}\int_\Omega \phi(\eta)d\eta\,,
\end{equation}
where 
$$
\alpha_n^{m+1} = \exp(-\sigma(x_n^{m+1}) \Delta t) = \exp(-\sigma(x_n^{m} + \Delta t v_n^m) \Delta t)\,, 
$$ 
while $\eta^{m+1}$ and $p^{m+1}$ are random variables in~\Cref{alg:f-RTE}.  

Consider a fixed particle index $n$. For different realizations of $\{\left( x_n^m\,,v_n^m \right) \}_{m=0}^{M}$, the value $\phi(v_n^m)$ changes with respect to the time index $m$. As a result, $\phi(v_n^m)$ can be treated as a discrete-time stochastic process on a probability space spanned by random variables $\{p^m,\eta^m\}_{m=1}^M$, and hence naturally generates a filtration. 

Moreover, we define the following conditional expectation for $m = 0,1,\ldots, M-1$
\[
    \tilde\bbE^m[~\cdot~|(x_n^m,v_n^m)] = \bbE_{\eta^\np} \bbE_{p^\np} \cdots \bbE_{\eta^M} \bbE_{p^M}[~\cdot~ | (x_n^m, v_n^m)] \,.
\]
Using these notations, the final objective functional~\eqref{eq:J2}, after forward Euler time discretization, can be approximated as sums of conditional expectations: 
\begin{equation}\label{eq:J_conditional}
    J \approx \tilde{J}^m:=  \frac{1}{N} \sum_{n=1}^N \mathbb{E}_{(x_n^{m},v_n^{m})\sim f(t_m,x,v ) } \left[ \tilde{\bbE}^{m} [r(x_n^M,v_n^M)  |(x_n^m,v_n^m)]  \right]\,,
\end{equation}
as a result of the law of total expectations. Note that the above formula holds for any $0\leq m \leq M-1$.

We can further write $\tilde{J}^m$ into $\tilde{J} ^m= \sum_{n=1}^N \tilde{J}^m_n$ following~\eqref{eq:J_conditional}, where
\begin{eqnarray}
\tilde{J}^m_n &=&  \frac{1}{N} \mathbb{E}_{(x_n^{m},v_n^{m})\sim f(t_m,x,v ) } \left[ \mathcal{R}_n^m\left(x_n^m, v_n^m\right)  \right] \, , \label{eq:Ji_conditional} \\
\mathcal{R}_n^m\left(x_n^m, v_n^m\right) &=&  \tilde{\bbE}^{m} \left[r\left(x_n^M,v_n^M\right) |(x_n^m,v_n^m) \right] \, . \label{eq:Ri_def}
\end{eqnarray}

The dependence of $\mathcal{R}_n^m\left(x_n^m, v_n^m\right)$ on the coefficient function $\sigma(x)$ is through the evaluations of $\{\sigma(x_n^{m+1})\}$ where $x_n^{m+1} = x_n^m + \Delta t\, v_n^m$, used in the acceptance-rejection probabilities in each $\mathbb{E}_{p^{m+1}}$; see~\eqref{eq:probab}.  Note that $\mathcal{R}_n^m$ is conditioned on $(x_n^m, v_n^m)$, so it can be seen as a function of $(x_n^m, v_n^m)$ and consequently a function of  $\sigma(x_n^{m+1})$ for a given coefficient $\sigma(x)$. When $(x_n^m, v_n^m)$ is considered a random variable, $\mathcal{R}_n^m$ is also a random variable.   Thus, using the \textbf{Score Function Method} for computing the Monte Carlo gradient, discussed in Section~\ref{subsubsec:score}, we can express the derivative of $\mathcal{R}_n^m$ with respect to each $ \theta_n^{m+1} := \sigma(x_n^{m+1})$ as
\begin{equation}\label{eq:dRidsigma}
    \frac{\partial \mathcal{R}_n^m}{\partial \theta_n^{m+1}} =  \tilde{\bbE}^{m} \left[ \frac{\partial \log \kappa_n^{m+1} }{\partial \theta_n^{m+1} }   r\left(x_n^M,v_n^M\right)  \right] =: G_n^m\, ,
\end{equation}
where the probability for the rejection sampling and its score function are
\[
\kappa_n^{m+1}  = 
    \begin{cases}
   \alpha_n^{m+1} ,  & \text{if\ }v_n^{m+1} = v_n^m  \\
    1- \alpha_n^{m+1} , & \text{otherwise}
    \end{cases}, \qquad    
 \frac{\partial \log \kappa_n^{m+1} }{\partial \theta_n^{m+1} }   = 
    \begin{cases}
   -\Delta t \,,  & \text{if\ }v_n^{m+1} = v_n^m  \\
   \Delta t \frac{ \alpha_n^{m+1}}{1-\alpha_n^{m+1}} \,, & \text{otherwise}
    \end{cases}.
\]
In~\eqref{eq:dRidsigma}, we fix a particular $\left(x_n^m,v_n^m\right)$. When $\left(x_n^m,v_n^m\right)$ is considered a random variable, both $G_n^m$ and $\mathcal{R}_n^m$ are also random variables.   Using the same trajectories from the Monte Carlo solver for the forward RTE (see~\Cref{alg:f-RTE}), we obtain samples $\{ \widehat {\mathcal{R}}_n^m\}_{m=0}^{M-1}$ according to~\eqref{eq:Ri_def}, and the Monte Carlo gradient $\{  \widehat {\mathcal{G}}_n^m\}_{m=0}^{M-1}$ based on~\eqref{eq:dRidsigma}, for  $\mathcal{R}_n^m$ and $G_n^m$, respectively:
\begin{equation}\label{eq:dRidsigma_discrete}
\widehat {\mathcal{R}}_n^m =   r\left(x_n^M,v_n^M\right),\qquad   \widehat {\mathcal{G}}_n^m  = 
     \begin{cases}
     -   r\left(x_n^M,v_n^M\right) \Delta t , & \text{if\ }  v_n^{m+1} = v_n^m , \\
      r\left(x_n^M,v_n^M\right) \dfrac{ \alpha_n^{m+1}  \Delta t }{1 - \alpha_n^{m+1} }, &\text{otherwise},
     \end{cases}
\end{equation}
with $m = 0,\ldots,M-1$.

Based on~\eqref{eq:Ji_conditional} and~\eqref{eq:Ri_def}, we can also treat $N^{-1} \widehat {\mathcal{R}}_n^m$ as a sample of $\tilde{J}^m_n$. Since all the derivations above are based on the same forward RTE particle trajectories, $\{(x_n^m, v_n^m)\}_{m=0}^M$, $n = 1,\ldots, N$, obtained by running~\Cref{alg:f-RTE}, the following holds for any $m = 0,\ldots, M-1$:
\begin{eqnarray*}
J & \stackrel{\text{time discretization}}{\approx} & \sum_{n=1}^N \tilde{J}^m_n \quad  \stackrel{\text{sample~\eqref{eq:Ji_conditional}}}{\approx} \quad \frac{1}{N}\sum_{n=1}^N  {\mathcal{R}}_n^m \quad  \stackrel{\text{sample~\eqref{eq:Ri_def}}}\approx \quad   \frac{1}{N}\sum_{n=1}^N \widehat {\mathcal{R}}_n^m,\\
  &\stackrel{\text{based on~\eqref{eq:J2_MC}}}\approx& \hat J  = \frac{1}{N}\sum_{n=1}^N r(x_n^M, v_n^M)\,.
\end{eqnarray*}
Thus, the Monte Carlo gradient  of  $\hat J$ with respect to $\sigma(x_n^{m+1})$ can be approximated by
\begin{equation}\label{eq:DTO_gradient_old}
\frac{\partial \hat J }{\partial  \theta_n^{m+1}  }\approx \frac{1}{N} \sum_{i=1}^N \frac{\partial  {\mathcal{R}}_{i}^m }{\partial  \theta_n^{m+1}  }  = \frac{1}{N} \frac{\partial \mathcal{R}_n^m }{\partial  \theta_n^{m+1}  } \approx  \frac{1}{N}\widehat{\mathcal{G}}_n^m,\quad \forall m\,,
\end{equation}
since $\mathcal{R}_{i}^{m}$ does not depend on $ \theta_n^{m+1} 
 = \sigma(x_n^{m+1})$ if $i\neq n$, $\forall m$. Combining \eqref{eq:DTO_gradient_old} with~\eqref{eq:dRidsigma_discrete}, we obtain a Monte Carlo gradient formula for $\frac{\partial \hat J }{\partial \theta_n^{m+1} }$.

However, it is worth noting that 
\begin{equation}\label{eq:two gradients}
 \frac{\partial \hat{J}}{\partial \theta_n^{m+1} }\neq \frac{\delta J}{\delta \sigma}(x_n^{m+1})\,, \quad \theta_n^{m+1} = \sigma(x_n^{m+1})\,.
\end{equation}
The left-hand side gradient treats $\theta_n^{m+1} = \sigma(x_n^{m+1})$ as a single parameter, and therefore the effective parameters are $MN$ scalars, $\{\theta_n^{m+1}\}$, where $n=1,\ldots,N$ and $m = 0,\ldots M-1$. On the other hand, $\frac{\delta J}{\delta \sigma}$ is a functional derivative with respect to the parameter function $\sigma(x)$ (the same with the OTD derivative~\eqref{eq:RTE_OTD_grad}). Since we assume $\sigma\in L^2$ and the dual space of $L^2$ is still $L^2$, the derivative $\frac{\delta J}{\delta \sigma}$ is also the gradient function $G(x)$. The right-hand side of~\eqref{eq:two gradients} is then function $G(x)$ evaluated at $x= x_n^{m+1}$. We use the following example to show how to relate both sides of~\eqref{eq:two gradients}.

Consider a separate mesh grid $\{\bar{x}_j\}$ in the spatial domain, and $Q_j$ is a small neighborhood of $\bar{x}_j$ for each $j$. Consider a perturbed parameter function $\sigma(x) + \delta \sigma(x)$, where $\delta \sigma(x) = \mathds{1}_{x \in Q_j}  \epsilon$ for some small constant $\epsilon$. We then have
\begin{equation}\label{eq:perturb1}
J(\sigma  +\delta \sigma ) - J(\sigma) \approx \int_x \frac{\delta J}{\delta \sigma}(x) \delta \sigma(x) \rd{x} = \epsilon  \int_{Q_j} \frac{\delta J}{\delta \sigma}(x) \rd{x}.
\end{equation}
On the other hand, with $\hat{J}(\sigma)$ denoting the value~\eqref{eq:J2_MC} calculated with a given parameter function $\sigma$, we have the following based on~\eqref{eq:DTO_gradient_old}:
\begin{equation}\label{eq:perturb2}
\hat{J}(\sigma  +\delta \sigma ) -  \hat{J}(\sigma) \approx  \sum_{n=1}^N \sum_{m=0}^{M-1}  \frac{\partial \hat{J} }{\partial \sigma(x_n^{m+1})} \delta \sigma(x_n^{m+1})  \approx 
\epsilon\,  \frac{1}{N} \sum_{n=1}^N \sum_{m=0}^{M-1}  \mathds{1}_{x_n^{m+1} \in Q_j} \,  \widehat{\mathcal{G}}_n^m\,,
\end{equation}
where we used~\eqref{eq:DTO_gradient_old}. Combining the last terms in~\eqref{eq:perturb1}-\eqref{eq:perturb2} and assuming the approximated gradient function $\frac{\delta J}{\delta \sigma}$ is piecewise constant on $\{Q_j\}$, we have 
\[
 \frac{\delta J}{\delta \sigma}(\bar x_j)\, |Q_j| \approx \int_{Q_j} \frac{\delta J}{\delta \sigma}(x) \rd{x} \approx \frac{1}{N} \sum_{n=1}^N \sum_{m=0}^{M-1}  \mathds{1}_{x_n^{m+1} \in Q_j} \, \widehat{\mathcal{G}}_i^n.
\]
Finally, we may approximate the gradient using values of~\eqref{eq:dRidsigma_discrete} by
\begin{eqnarray}
 \frac{\delta J}{\delta \sigma}(\bar x_j)  & \approx& \frac{1}{|Q_j|}  \frac{1}{N}  \sum_{n=1}^N \sum_{m=0}^{M-1}  \mathds{1}_{x_n^{m+1} \in Q_j}   \widehat{\mathcal{G}}_n^m\nonumber
 \\
 &=&\frac{1}{|Q_j|} \frac{ \Delta t}{N}\sum_{n=1}^N\sum_{m=1}^{M} \ \mathds{1}_{x_n^{m} \in Q_j} \  r(x_n^M, v_n^M) \ \xi_n^m,\quad \xi_n^m = \begin{cases}
 -1 , & \text{if\ }  v_n^{m} = v_n^{m-1}, \\
 \frac{\alpha_n^m}{1 - \alpha_n^m } , &\text{otherwise}.
 \end{cases} \label{eq:DTO_grad}
\end{eqnarray}
Note that sorting is needed in computing \eqref{eq:DTO_grad} for all particles in the long-time horizon.

In this DTO approach, the adjoint particle does  not show up explicitly, but is implicitly represented by the particle representation of the gradient, $\widehat{\mathcal{G}}_n^m$, where $1\leq n\leq N$, $0\leq m \leq M-1$. Started from the DTO approach, this is our second example of using an adjoint Monte Carlo method for efficient gradient computation.

\subsection{Numerical Examples}
Next, we show some numerical examples to demonstrate the adjoint Monte Carlo method for computing the gradient of the objective functional with respect to the parameter. We will also show the result based on the finite volume method (FVM) for comparison.

Consider a 1D spatial domain $[-2,2]$ and the velocity domain on a unit interval $[-1,1]$. The time domain is $[0,0.5]$. The spatial domain is uniformly discretized with the spacing $\Delta x = 0.05$ at which the gradient of the RTE-constrained optimization problem with respect to the scattering coefficient $\sigma(x)$ is evaluated. The objective function is in the form of~\eqref{eq:J2} with $r(x,v) = |v|^2 \mathbbm{1}_{x<0}$. The initial distribution $f_0(x,v) = \frac{1}{\sqrt{\pi}}\left(\exp(-4(x-0.5)^2) + \exp(-4(x+0.5)^2) \right)$.   The parameter where the gradient is evaluated is $\sigma(x) = 2+2\exp(-4x^2)$. 

\begin{figure}
\centering
\subfloat[The approximated gradient]{\includegraphics[width = 0.49\textwidth]{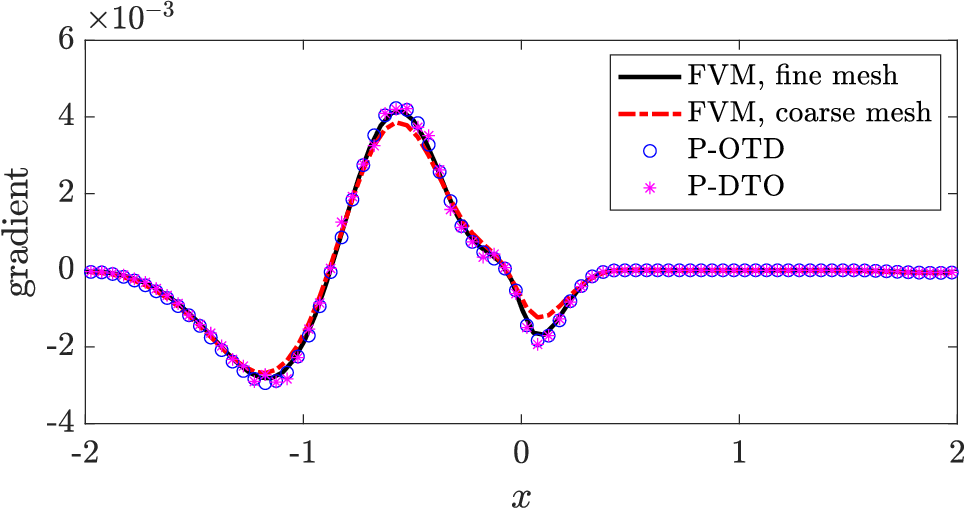}\label{fig:RTE_grad}}
\subfloat[The standard deviation of adjoint MC methods]{\includegraphics[width = 0.49\textwidth]{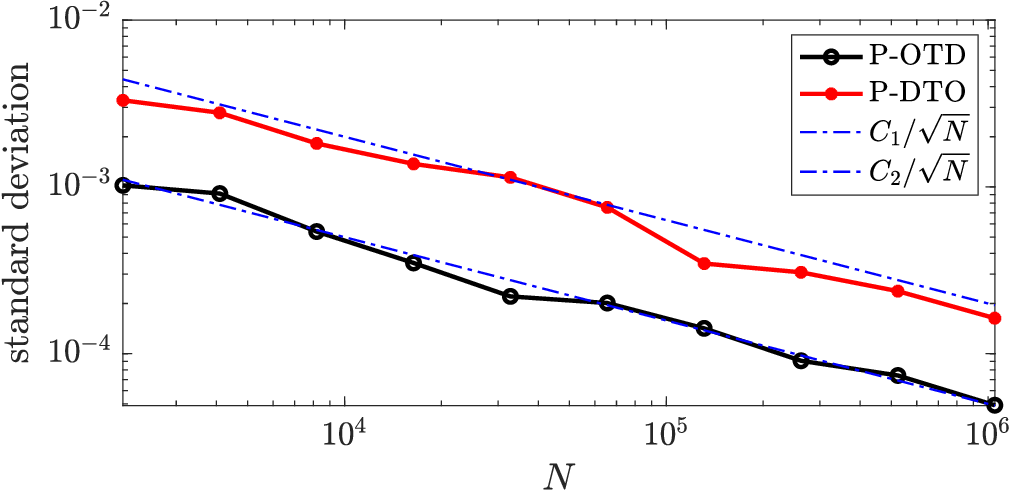}\label{fig:RTE_std}}
\caption{Illustration of the approximated gradient for an RTE-constrained optimization problem based on various numerical methods (left) and the standard deviation for the adjoint Monte Carlo methods based on the OTD and DTO approaches (right).}
\end{figure}

For FVM, we discretize the spatial and time domain with $\Delta x = 0.05$, $\Delta t = 0.01$ (coarse) and  $\Delta x = 0.005$, $\Delta t = 0.001$ (fine) respectively. In the fine-grid discretization, we have to reduce the time spacing due to the CFL condition and the use of the forward Euler time integrator.  The velocity domain $[-1, 1]$ is uniformly discretized based on the spacing $\Delta v = 0.05$. For the adjoint Monte Carlo method based on the OTD approach (denoted as ``P-OTD''),  we use $N = 10^6$ particles and time spacing $\Delta t = 0.01$, following the procedure detailed in Section~\ref{subsec:P-OTD}. 
For the adjoint Monte Carlo method based on the DTO approach (denoted as ``P-DTO''),  we also use $N = 10^6$ particles and time spacing $\Delta t = 0.01$, following the scheme presented in Section~\ref{subsec:RTE_DTO}.  We average the results of adjoint Monte Carlo methods from $64$ i.i.d.~runs to illustrate the final gradient.   We show the approximated gradients from these four methods in Figure~\ref{fig:RTE_grad}.  The values obtained by the adjoint Monte Carlo methods coincide with the ones from FVM with a fine discretization. However, the result obtained from the FVM with a coarse discretization is subject to visible numerical error.

In Figure~\ref{fig:RTE_std}, we range the number of particles in the adjoint Monte Carlo methods from $N = 2^{11}$ to $2^{20}$ and then compute the standard deviation of the  $L^2$ error of the gradient based on the $64$ i.i.d.~runs.  The gradient computed using the FVM with a fine discretization is regarded as the ground truth.   Both the P-OTD and P-DTO methods are subject to the standard Monte Carlo error $\mathcal{O}(N^{-1/2})$.  It is worth noting that the P-DTO method is subject to a larger random error than the P-OTD method for the RTE-constrained optimization.

\section{Case Study II: Boltzmann Equation}\label{sec:Boltzmann}
In this section, we present an adjoint Monte Carlo method for the Boltzmann equation, a famous kinetic equation with a quadratic nonlinearity. The derivation here are more intricate as we need to incorporate the binary collision feature of the forward problem into the adjoint particle system.

For simplicity, we consider the spatially homogeneous Boltzmann equation
\begin{equation}\label{eq:homoBoltz}
    \frac{\partial f}{\partial t}  =  Q(f,f)\, .
\end{equation}
The nonlinear collision operator $Q(f,f)$, which describes the binary collisions among particles, is defined as
\begin{equation*}%
    Q(f,f) = \int_{\bbR^3} \int_{\cS^2} q(v-v_1,\sigma) (f(v_1')f(v') - f(v_1) f(v)) \rd \sigma \rd v_1,
\end{equation*}
in which  $(v',v_1')$ represent the post-collisional velocities associated with the pre-collisional velocities $(v,v_1)$, the collision kernel function $q\geq 0$, and the $\sigma$ integral is over the surface of the unit sphere $\cS^2$.

By conserving the momentum $v+v_1$ and the energy $v^2 + v_1^2$, we have
\begin{equation}\label{BoltzmannSolution}
\begin{split}
v' &= 1/2 (v+v_1) + 1/2 |v-v_1| \sigma,  \\
v_1'&=1/2 (v+v_1) - 1/2  |v-v_1| \sigma,
\end{split}
\end{equation}
where $\sigma \in \cS^2$ is a collision parameter. We will hereafter use the shorthand notation 
\begin{eqnarray*}%
 (f, f_1, f', f_1') &=&  (  f(v), f(v_1), f(v'), f(v_1') ), \\
 (\hat f, \hat f_1,\hat f',  \hat f_1') &=&  \rho^{-1}(  f(v), f(v_1), f(v'), f(v_1') )
\end{eqnarray*}
for the values of $f$ in the Boltzmann collision operator $Q$, as well as the normalized densities $\hat f$ with $\rho = \int f(v) \rd v$ and $\int \hat f(v) dv=1$.

Physical symmetries imply that
\begin{equation}~\label{eq:GeneralKernel}
    q(v-v_1,\sigma) =  \tilde q(|v-v_1|,\theta),
\end{equation}
where $\cos \theta = \sigma\cdot \alpha$ with $\alpha = \frac{v-v_1}{|v-v_1|}$ and $\sigma = \frac{v'-v'_1}{|v'-v'_1|} = \frac{v'-v'_1}{|v-v_1|}$. That is, the collision kernel we consider here only depends on the size of the pre-collision relative velocity $|v-v_1|$, and the angle $\theta$ between the pre- and post-relative velocity directions.

In~\cite{caflisch2021adjoint}, the collision rule~\eqref{BoltzmannSolution} was re-written conveniently in terms of matrix-vector multiplications, as
\begin{equation}\label{eq:AB_vel_General}
        \begin{pmatrix}
    v'\\
   v_1'
    \end{pmatrix}    = A(\sigma,\alpha) \begin{pmatrix}
    v\\
    v_1
     \end{pmatrix},
   \quad 
\end{equation}
where
\begin{equation} \label{eq:AB_def}
A(\sigma ,\alpha )= \dfrac{1}{2}\begin{pmatrix}
    I+\sigma  \alpha ^T & I-\sigma  \alpha ^T\\
    I-\sigma  \alpha ^T & I+\sigma  \alpha ^T
    \end{pmatrix},\quad 
\end{equation}
in which $I$ is the identity matrix in $\Rthree$. It is worth noting that 
$$
A(\sigma ,\alpha )^\top = A(\sigma ,\alpha )^{-1}= : B(\sigma ,\alpha ).
$$
We highlight the dependence of $A$ and $B$ on $\sigma$ and $\alpha$ by writing $A(\sigma ,\alpha )$ and $B(\sigma ,\alpha )$. Here, for a fixed pair of velocity $(v,v_1)$, we also define the Jacobian matrix $C\in\mathbb{R}^{6\times 6}$ given by
\begin{equation}\label{eq:C_def}
C(\sigma,\alpha) = \dfrac{\partial (v',v_1')}{\partial (v,v_1)}.
\end{equation}
Even though~\eqref{eq:AB_vel_General} is a matrix-vector multiplication, $C$ and $A$ do not always coincide since the collision parameter $\sigma$ can depend on the velocity pair $(v,v_1)$ for certain collision kernel $q$.

Next, we will assume that the collision kernel $q$ is uniformly bounded from above. For some $\Sigma>0$, we have
\begin{equation}~\label{eq:UpperBound}
    q(v-v_1,\sigma) \leq \Sigma
\end{equation}
for all  $v$, $v_1$, and $\sigma$. In practice, this is not a restriction because one can take $\Sigma$ to be the largest value of $q$ for the discrete set of velocity and parameter values.

Although the formulation and analysis presented here are valid for the general collision model~\eqref{eq:GeneralKernel}, the following form is often considered,
\begin{equation}\label{eq:VSS/VHS}
    q(v-v_1,\sigma) =  \tilde q(|v-v_1|,\theta) = C_\kappa(\theta) |v-v_1|^\beta\,.
\end{equation}
We refer to~\cite[Sec.~1.4]{villani2002review} for more modeling intuition regarding this type of collision models. The collision model~\eqref{eq:VSS/VHS} accommodates both the variable hard sphere (VHS) model, when $C_\kappa(\theta)$ is constant, and the variable soft sphere (VSS) model~\cite{KouraMatsumoto} with angular dependence.

With the assumption~\eqref{eq:UpperBound}, the Boltzmann equation~\eqref{eq:homoBoltz} can be further written as
\begin{eqnarray}
    \frac{\partial f}{\partial t}   &= & \int_{\Rthree}\int_\Stwo \left(f'f_1' - f f_1\right)  q (v-v_1,\sigma) \rd \sigma \rd v_1 \nonumber \\
    & =& \iint f'f_1' q  \rd \sigma \rd v_1 - \iint ff_1 q  \rd\sigma \rd v_1\nonumber  \\
    & = & \iint f'f_1' q  \rd \sigma \rd v_1  + \iint ff_1 (\Sigma  - q)  \rd\sigma \rd v_1  - \iint \Sigma  ff_1 \rd\sigma \rd v_1 \nonumber\\
    &=& \iint f'f_1' q  \rd \sigma \rd v_1  + f \iint f_1 (\Sigma - q)  \rd \sigma \rd v_1  - \mu f. 
 \label{eq:VHS1}
\end{eqnarray}
Recall $\rho = \int f_1 \rd v_1$ is the total mass. We define ${\bar A} = \int_{0}^{2\pi}\int_{0}^{\pi} \sin(\theta) \rd \theta \rd \varphi = 4\pi$, the surface area of the unit sphere, and then we have
\begin{equation} \label{eq:mu}
\mu = {\bar A} \Sigma \rho.
\end{equation}
If we multiply both sides of~\eqref{eq:VHS1} by an arbitrary test function $\phi(v)$, then divide by $\rho$ and finally integrate over $v$, we obtain 
\begin{eqnarray*}
       \frac{\partial ( \int \phi \hat f \rd v)}{\partial t}  &= & \rho \iiint  \phi \hat f' \hat f_1' q  \rd \sigma \rd v_1 dv + \rho \iiint \phi \hat f \hat f_1 (\Sigma  - q)  \rd \sigma \rd  v_1 \rd v -  \mu \int  \phi \hat f \rd  v\\
       & = &\rho \iiint  \phi' \hat f \hat f_1 q  \rd  \sigma  \rd 
 v_1 dv +\rho \iiint  \phi \hat f \hat f_1 (\Sigma  - q)  \rd  \sigma \rd  v_1 \rd  v -  \mu \int  \phi \hat f \rd v,
\end{eqnarray*}
using $\rd v \rd v_1=\rd v' \rd  v_1' $ and then interchanging notation $(v,v_1)$ and $(v',v_1')$ in the first integral.
After applying the explicit Euler time integration from $t_k$ to $t_{k+1} = t_k + \Delta t$, we obtain
\begin{eqnarray}
     && \int \phi\ \hat f (v,t_{k+1}) dv \nonumber \\
       &=&   (1-\Delta t \mu ) \int \phi\ \hat f(v,t_{k}) \rd v  \label{eq:non_collide} \\
        &&  + \iiint  \phi'  \ {\Delta t \mu \Sigma^{-1}  {q}} \  \hat f(v,t_{k})  \hat f_1(v_1,t_{k})   {\bar A}^{-1} \rd 
        \sigma \rd  v_1 \rd v  \label{eq:collide} \\
       &&  +  \iiint \phi\ {\Delta t  \mu \Sigma^{-1} (\Sigma -q)} \  \hat f(v,t_{k})  \hat f_1(v_1,t_{k})  {\bar A}^{-1} \rd 
        \sigma \rd  v_1 \rd v. \label{eq:v_collide}
\end{eqnarray}
For a fixed $k$, $\hat f(v,t_k) =  f(v,t_k)  / \rho$ is the probability density of $v$, and $ \hat f(v,t_k)  \hat f_1(v_1,t_k) {\bar A}^{-1}$ is the probability density over the three variables $(v,v_1,\sigma)$. The three terms~\eqref{eq:non_collide}, \eqref{eq:collide} and~\eqref{eq:v_collide} on the right-hand side of this equation represent the sampling of the collision process, using rejection sampling, as described in the DSMC algorithm in~\Cref{subsec:DSMC}.

\subsection{The (Forward) DSMC Method}\label{subsec:DSMC}

\begin{algorithm}
  \caption{DSMC Algorithm Using Rejection Sampling for a General Collision Kernel $q$ in the Nonlinear Boltzmann Equation.}\label{alg:VHS_DSMC} 
\begin{algorithmic}[1]
\State  Compute the initial velocity particles based on the initial condition, $\mathcal V^0 = \{v^0_1,\dots,v^0_N\}$. 
\For{$k=0$ to $M-1$}
\State Given $\mathcal V^{k}$, choose $N_c=\ceil[\big]{ \Delta t \mu N/2}$  velocity pairs $(i_\ell ,i_{\ell_1})$ uniformly without replacement. The remaining $N-2N_c$ particles do not have a virtual (or real) collision and set $v^{k+1}_{i} = v^{k}_{i}$.
\For{$\ell = 1$ to $N_c$}
\State Sample $\sigma^k_\ell$ uniformly over $\Stwo$.
\State Compute $\theta^k_\ell = \arccos ( \sigma^k_\ell \cdot  \alpha^k_\ell )$ and $q^k_\ell=q( |v^k_{i_{\ell }}-v^k_{i_{\ell_1}}|, \theta^k_\ell)$.
\State Draw a random number $\xi^k_\ell$ from the uniform distribution $\mathcal{U}([0,1])$.
\If{ $\xi^k_\ell \leq q^k_\ell /\Sigma $}
\State Perform real collision for $( v^k_{i_{\ell }}, v^k_{i_{\ell_1}})$ following~\eqref{BoltzmannSolution} and obtain $({v^k_{i_{\ell }}}', {v^k_{i_{\ell_1}}}')$.
\State Set $( v^{k+1}_{i_{\ell }} , v^{k+1}_{i_{\ell_1 }} ) = ({v^k_{i_{\ell }}}', {v^k_{i_{\ell_1}}}'$).
\Else 
\State The virtual collision is not a real collision. Set $( v^{k+1}_{i_{\ell }} , v^{k+1}_{i_{\ell_1}} ) = ({v^k_{i_{\ell }}}, {v^k_{i_{\ell_1}}}$).
\EndIf
\EndFor
\EndFor
\end{algorithmic}
\end{algorithm}

The study of rarefied gas dynamics, particularly in regimes where the continuum assumption is not valid, has always posed significant challenges. Over the past few decades, the Direct Simulation Monte Carlo (DSMC) method has emerged as a pivotal tool in this domain, offering a particle-based, probabilistic approach to tackle the complexities of the Boltzmann equation.

At the heart of DSMC lies the representation of gas molecules as simulated particles. Each of these particles embodies a multitude of real molecules, moving and interacting based on probabilistic models. The method decouples molecular motion and collisions, simplifying the simulation process. A cornerstone of the DSMC method is its treatment of molecular collisions. Through probabilistic models, DSMC captures the essence of intermolecular interactions, ensuring that the simulated dynamics align closely with real-world physics. 
The inherent parallel nature of DSMC makes it well-suited for high-performance computing environments.  The versatility of DSMC has led to its adoption in a wide range of applications from aerodynamics to  microfluidics, including providing the numerical solution to the Boltzmann equation~\eqref{eq:homoBoltz}. 

Bird's book~\cite{bird1994molecular} is often considered the seminal reference for the DSMC method. A comprehensive review of the DSMC method and its applications can also be found 
in~\cite{alexander1997direct}. The DSMC tutorial written by Pareschi and Russo~\cite{pareschi2001introduction} is a very clear manual for the application of DSMC for the Boltzmann equation. The book~\cite{cercignani2013mathematical} by Cercignani, Illner  and Pulvirenti provides a mathematical foundation for the Boltzmann equation and methods like DSMC. For the rest of the section, we focus on the DSMC method particularly designed for the Boltzmann equation. Instead of directly solving the Boltzmann equation, DSMC simulates the motion and collisions of a finite number of representative gas molecules.

In the DSMC method~\cite{bird1970direct,nanbu1980direct,babovsky1989convergence}, we consider a set of $N$ velocities evolving in discrete time due to collisions whose distribution can be described by the distribution function $f$ in~\eqref{eq:homoBoltz}. We divide time interval $[0,T]$ into $M$ number of sub-intervals of size $\Delta t=T/M$. At the $k$-th time interval, the particle velocities are represented as
\begin{equation*}%
{\cal V}^k =\{v_1, \ldots , v_N\} (t_k),
\end{equation*}
and we denote the $i$-th velocity particle in ${\cal V}^k $ as $v_i^k$. The distribution function $f(v,t)$ is then discretized by the empirical distribution
\begin{equation}\label{eq:empirical}
    f(v,t_k) \approx \frac{\rho}{N} \sum_{i=1}^N \delta(v-v_i^k), \quad k = 0,\ldots,M.
\end{equation}

We define  the total number of virtual collision pairs $N_c =\ceil[\big]{ \Delta t \mu N/2}$. Note that the number of particles having a virtual collision is $2N_c \approx \Delta t \mu N$. Thus, the probability of having a virtual collision is $\Delta t \mu$, and the probability of not having a virtual collision is $1-\Delta t \mu$,  matching the term~\eqref{eq:non_collide}.  
For each velocity $v_i^k\in {\cal V}^k $ in this algorithm,  there are three possible outcomes, whose probabilities are denoted by $h_j$ where $j=1,2,3$:
\begin{center}
\begin{tabular}{ll}
\textbf{Outcome} &   \textbf{Probability} \\
1. No virtual collision  & $h_1 = 1-\Delta t \mu$  \\
2. A real collision      & $h_2 = \Delta t \mu \  q_i^k /\Sigma  $ \\
3. A virtual, but not a real, collision  &  $h_3 = \Delta t \mu (1- q_i^k /\Sigma)$ \\
\end{tabular}
\end{center}
Here, $q_i^k = q(v_i^k-v_{i_1}^k,\sigma_i^k)$ where $v_{i_1}^k$ represents the virtual collision partner of $v_i^k$ and $\sigma_i^k$ is the sampled collision parameter for this pair. Note that the total probability of no real collision is $h_1+h_3= 1 - \Delta t \mu q_i^k/\Sigma$. Since the $h_j$'s depend on the collision kernel $q$, we may also view them each as function of $v-v_1$ and $\sigma$. Later, we will use the following score functions based on the definition in~\eqref{eq:score_def}:
\begin{eqnarray*}
\partial_{v_i^k} (\log h_1) &=& 0, \nonumber  \\
\partial_{v_i^k} (\log h_2) &=& (q_i^k)^{-1}  \partial_{v_i^k} q_i^k, \\
\partial_{v_i^k} (\log h_3)&=& - (\Sigma -q_i^k)^{-1}  \partial_{v_i^k} q_i^k \,. \nonumber
\end{eqnarray*}
We will combine the Score Function Method (Section~\ref{subsubsec:score}) with the adjoint-state method (Section~\ref{subsec:adjoint}) to design the adjoint Monte Carlo algorithm. Recall $\rho$, $\Delta t$ and $\Sigma$ are constants. Also note that $\partial_{v_{i_1}^k} (\log h_j)= - \partial_{v_{i}^k} (\log h_j)$, $j=1,2,3$, if $(v_{i}^k,v_{i_1}^k)$ is a virtual collision pair.

We can decide whether a particle participates in a \textit{virtual} collision or not through uniform sampling. However, in order to determine whether a selected virtual collision pair participates in the  \textit{actual} collision or not, we need to use the rejection sampling since the probability $q_i^k$ is pair-dependent. We further remark that if a virtual velocity pair is rejected for a real collision (Outcome \#2), it is automatically accepted for a virtual but not real collision (Outcome \#3). With all the notations defined above, we present the DSMC algorithm using the rejection sampling in~\Cref{alg:VHS_DSMC}. The algorithm applies to any general collision kernel satisfying~\eqref{eq:UpperBound}. 

Note that the DSMC algorithm is a single sample of the dynamics for $N$ particles following equations \eqref{eq:non_collide}-\eqref{eq:v_collide}. It is not the same as $N$ samples of a single particle because of the nonlinearity of these equations. For consistency, the derivation of the adjoint equations will also employ a single sample of the $N$-particle dynamics.

\subsection{DTO Approach: The Adjoint DSMC Method }\label{subsec:adjoint_DSMC}
As a test problem, we consider an optimization problem for the spatially homogeneous Boltzmann equation~\eqref{eq:homoBoltz}. 
The initial condition is
\begin{equation}
f(v,0)=f_0(v;\theta), \label{BoltzmannIC}
\end{equation}
where $f_0$ is the prescribed initial data depending on the parameter $\theta$, which can be high-dimensional. The goal is to find $\theta$ which optimizes the objective function
\begin{equation}\label{eq:OTD_obj}
J (\theta) = \int_{\bbR^3} \phi (v) f(v,T) \rd v,
\end{equation}
where $f(v,T)$ is the solution at time $T$ to~\eqref{eq:homoBoltz} given the initial condition~\eqref{BoltzmannIC}, and thus $f(v,T)$ depends on $\theta$ through the initial condition. This is a PDE-constrained optimization problem. The adjoint DSMC method proposed in~\cite{caflisch2021adjoint} is an efficient particle-based method to compute the gradient of the objective function~\eqref{eq:OTD_obj} based on the forward DSMC scheme (\Cref{subsec:DSMC}). It was generalized to a much more general collision kernel in~\cite{yang2023adjoint}. We will briefly review its derivation and the main result.

We first rewrite the objective function~\eqref{eq:OTD_obj} at time $T= t_M$ as
\begin{equation} \label{eq:obj1}
J = \mathbb{E}\left[ \bar{\phi}^M \right],\quad \text{where} \quad 
\bar{\phi}^M = \frac{\rho}{N}\sum_{i=1}^{N} \phi_i^M,
\end{equation}
where $\rho = \int f(v,T) \rd v$, $\phi_i^M=\phi(v_i^M)$ with $v_i^M \sim f(v,T)$. The expectation $\mathbb{E}$ is taken over all the randomness in the forward DSMC simulation over $M$ time steps.

Next, we define the expectations for each step of the forward DSMC algorithm. For simplicity, we assume that the number of particles $N$ is even. The velocities change at a discrete time in the DSMC~\Cref{alg:VHS_DSMC}. This involves the following two steps at time $t_k$ where $k=0,\ldots,M-1$.

\begin{itemize}
    \item[Step 1. ] Randomly (and uniformly) select collision pairs, $v_i^k$ and $v_{i_1}^k$, and also to randomly select collision parameters. The expectation over this step will be denoted as $\mathbb{E}_p^k$ (with ``p" signifying collision ``parameters" and collision ``pairs"). 
    \item[Step 2. ] Perform collisions at the correct rate, i.e., the given collision kernel, $q(v -v_1, \sigma)$, using the rejection sampling. This is performed by choosing outcome $j$ with probability \begin{equation}\label{eq:rejection_hij}
        h^k_{ij}=h_j (v^k_i -v^k_{i_1}, \sigma_i^k)
    \end{equation} 
    for $j=1,2,3$ regarding this collision pair $(v_i^k, v_{i_1}^k)$; see~\Cref{subsec:DSMC} for the definition of $\{h_j\}$. The expectation for this step will be denoted as $\mathbb{E}_r^k$ (with ``r" signifying collision ``rejection sampling"). 
\end{itemize}
The total expectation over the step from time $t_k$ to $t_{k+1}$ is
\begin{equation*}%
\mathbb{E}^k = \mathbb{E}_p^k \, \mathbb{E}_r^k  .
\end{equation*}

Two ways of deriving the same adjoint DSMC method were presented in~\cite{yang2023adjoint}. We will  describe below the Lagrangian approach following the general framework in~\Cref{subsec:adjoint}.
Consider the Lagrangian
\begin{eqnarray}
 L\left(\{v_i^k\},\{\gamma_i^k\},\theta \right) &= &  J + \frac{1}{2}\sum_{k=0}^{M-1} \sum_{i=1}^N \mathbb{E}^k \left[   \begin{pmatrix} \gamma_i^{k+1} \\ \gamma_{i_1}^{k+1} \end{pmatrix} \cdot  \left( \cC_{i}^k \begin{pmatrix} v_i^k \\ v_{i_1}^k \end{pmatrix} - \begin{pmatrix} v_i^{k+1} \\ v_{i_1}^{k+1} \end{pmatrix} \right)   \bigg| \begin{pmatrix} v_i^k \\ v_{i_1}^k \end{pmatrix} \right] +  \nonumber \\
&&   \sum_{i=1}^N \mathbb{E}_{\hat{v}_{0i} \sim f_0(v;\theta)} \left[  \gamma_{i}^0 \cdot  \left( \hat{v}_{0i} - v_i^0 \right) \right]. \label{eq:full_Lag}
\end{eqnarray}
The ``$1/2$'' scaling is to avoid enforcing the collision rule twice. Note that $v_i^{k}$ and $v_{i_1}^{k}$ are a collision pair, but that neither $(v_i^{k+1}, v_{i_1}^{k+1})$ nor $(v_i^{M}, v_{i_1}^{M})$  are a collision pair. In this application of the adjoint-state method (see Section~\ref{subsec:adjoint}), the state variables are the random velocity particles $\{v_i^{k}\} \subset \bbR^3$, and the adjoint variables are $\{\gamma_i^{k}\}\subset \bbR^3$, which are also random. 

In the second step, we first select a time step $t_k$, condition the Lagrangian $L$ on $\mathcal{V}^k = \{v_i^{k}\}_{i=1}^N$, and  differentiate it
with respect to each $v_i^k$, $i=1,\ldots, N$, to obtain the corresponding adjoint equations. For a given $k = 0,\ldots, M-1$, and a particular velocity $v_i^{k}$, the only elements in ${\cal V}^{k+1}$ that depend on $v_i^k$ are $v_i^{k+1}$ and $v_{i_1}^{k+1}$, which are determined by 
\begin{equation} \label{eq:general_collision}
\begin{pmatrix}
v_i^{k+1}\\
v_{i_1}^{k+1}
\end{pmatrix} = 
\mathcal{C}^k_{i} 
\begin{pmatrix}
v_i^{k}\\
v_{i_1}^{k}
\end{pmatrix},
\end{equation}
where the operator $\mathcal{C}_{i}^k$ is one of the following three possible matrices each with probability $h_{ij}^k$ given in~\eqref{eq:rejection_hij} for $j=1,2,3$,
\begin{equation}\label{eq:collision_operator}
{\cal C}_i^k = 
    \begin{cases}
    I,  & j=1, \text{$(v_i^k,v_{i_1}^k)$ does not have a real or virtual collision}, \\
     A(\sigma_i^k, \alpha_i^k), & j=2,  \text{$(v_i^k,v_{i_1}^k)$ has a real collision},\\
    I, & j=3,  \text{$(v_i^k,v_{i_1}^k)$ has a virtual, but not a real, collision},
    \end{cases}
\end{equation}
where $I \in \mathbb{R}^{6\times 6}$ is the identity matrix.  Note that $A(\sigma,\alpha)$ and $\{h_{ij}^k\}$ are defined in~\eqref{eq:AB_def} and~\eqref{eq:rejection_hij}, respectively. Furthermore, we find the adjoint of the Jacobian matrices $D$ given by
\begin{equation} \label{eq:back_D}
    D =  \left[ \dfrac{\partial (v',v'_{1})}{\partial (v,v_1)} \right]^\top = \begin{cases}
    I,  & j=1, \\
    [C(\sigma, \alpha)]^\top , & j=2,\\
    I, & j=3\,,
    \end{cases}
\end{equation}
where $C(\sigma,\alpha)$ is defined in~\eqref{eq:C_def}, and is not necessarily the same as $A(\sigma,\alpha)$, especially when the collision kernel is angle dependent.

When we differentiate the Lagrangian $L$ with respect to $v_i^k$ and when $k\neq 0$ or $M-1$, the only term that depends on $v_i^k$ is
\[
\mathbb{E}^k \left[   \begin{pmatrix} \gamma_i^{k+1} \\ \gamma_{i_1}^{k+1} \end{pmatrix} \cdot  \left( \cC_{i}^k \begin{pmatrix} v_i^k \\ v_{i_1}^k \end{pmatrix} - \begin{pmatrix} v_i^{k+1} \\ v_{i_1}^{k+1} \end{pmatrix} \right)   \bigg| \begin{pmatrix} v_i^k \\ v_{i_1}^k \end{pmatrix} \right] \,.
\]
It is crucial to note that there are two places in this term that potentially depend on $v_i^k$: (1) the $v_i^k$ inside ``$[\,\cdot\,]$'' and (2) the probability $h_{ij}^k$ given in~\eqref{eq:rejection_hij} used to define $\mathbb{E}_r^k$, which is part of $\mathbb{E}^k$. Once we capture these two contributions, it follows that, for any test function $\varphi$,
\begin{equation}\label{eq:phi_diff_final}
    \frac{\partial}{\partial v_i^k}  \mathbb{E}^k [ \varphi  ]   =  \frac{\partial}{\partial v_i^k}  \mathbb{E}_{p}^k \mathbb{E}_{r}^k [ \varphi  ]  = 
    \mathbb{E}^k \left[ \frac{\partial}{\partial v_i^k}  \varphi + \varphi \left(\frac{\partial}{\partial v_i^k} \log h_i^k \right) \right]\,,
\end{equation}
where $h_i^k = h_{ij}^k$ with probability $h_{ij}^k$. We address that, similar to~\eqref{eq:score_grad1}, the score function $\frac{\partial}{\partial v_i^k} \log h_i^k$ shows up. Here, we have applied techniques in the \textbf{score function Monte Carlo estimator} introduced earlier in~\Cref{subsubsec:score}.

\begin{figure}
\begin{center}
\begin{tikzpicture}
\centering
\draw[black, thick,dashed,-latex] (1,.2) -- (3,2);
\draw[black, thick,-latex] (3,2) -- (7,4);
\draw[blue, thick,-latex]   (6.8,3.8)-- (3.1,1.95);
\draw[black, thick,-latex] (5,.2) --(3,2) ;
\draw[blue, thick,-latex] (3.15, 2) --(5.15,0.2) ;
\draw[black, dashed, thick,-latex] (3,2) -- (2,3);
\draw[black, thick,-latex] (7,4) -- (10,7);
\draw[blue, thick,-latex] (9.85,7) -- (7,4.2) ;
\draw[black, thick,-latex] (7,4) -- (4,7);
\draw[blue, thick,-latex] (4.15,7) -- (6.95,4.2);
\draw[black, thick,-latex] (10,2) -- (7,4) ;
\draw[blue, thick,-latex] (7.15,3.8) -- (9.85,2);
\draw[black, thick,-latex] (8,.2) -- (10,2);
\draw[blue, thick,-latex] (9.85,2) -- (7.85,.2);
\draw[black, thick,dashed,-latex] (12,.2) -- (10,2) ;
\draw[black, dashed, thick,-latex] (10,2) -- (11.5,2.8);
\node(x) [gray]  at (0,1.4) {$t_{k-1}$};
\node(y) [gray]  at (0,4) {$t_{k}$};
\node(z) [gray]  at (0,7) {$t_{k+1}$};
\draw[gray, thick,-stealth] (x) -- (y);
\draw[gray, thick,-stealth] (y) -- (z);
\node (a) at (3.3,1) {$v_i^{k-1}$};
\node [blue] at (5,1) {$\gamma_i^{k-1}$};
\node (a) at (4,3.1) {$v^{k}_i$};
\node [blue]  at (6,3.1) {$\gamma_i^{k}$};
\node (a) at (4,6) {$v_i^{k+1}$};
\node [blue]  at (6,6) {$\gamma_i^{k+1}$};
\node [blue]  (a) at (8.2,1) {$\gamma_{i_1}^{k-1}$};
\node (a) at (9.5,1) {$v_{i_1}^{k-1}$};
\node (a) at (9,3.1) {${v}_{i_1}^{k}$};
\node [blue]  at (7.5,3.1) {$\gamma_{i_1}^k$};
\node (a) at (9.6,6) {$v_{i_1}^{k+1}$};
\node [blue]  (a) at (8.2,6) {$\gamma_{i_1}^{k+1}$};
\end{tikzpicture}
\caption{An illustration for the Boltzmann forward and adjoint Monte Carlo methods at intermediate time steps $t_{k-1}$, $t_k$ and $t_{k+1}$: the forward velocity particles $\{v_i^k\}$ (black) evolve forward in time while the adjoint particles $\{\gamma_i^k\}$ (blue) propagate backward in time, both along the same trajectories sharing the randomness from the forward DSMC sampling.\label{fig:DSMC}}
\end{center}
\end{figure}
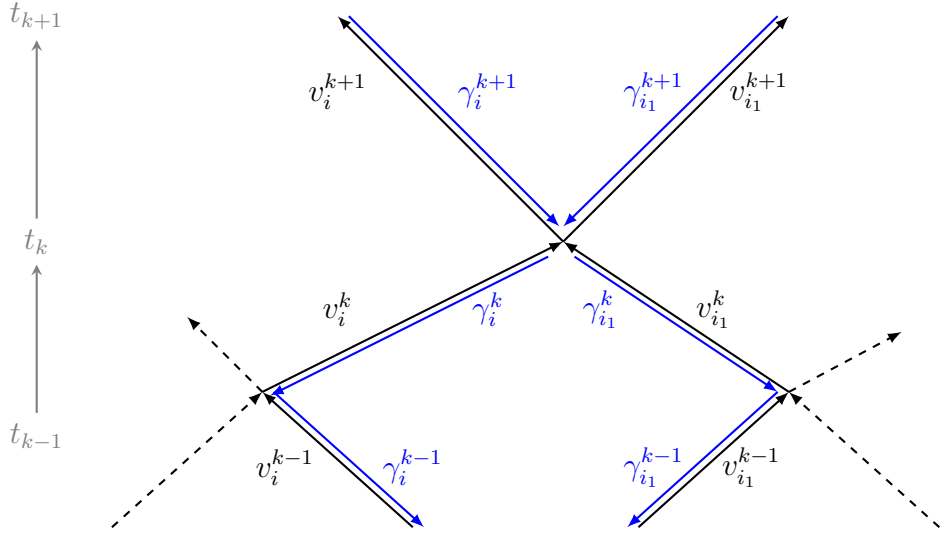

The introduction of the score function is crucial and the other steps are relatively straightforward. Thus, we omit the remaining derivation details and refer interested readers to~\cite{yang2023adjoint}. We state the final adjoint equations with the score functions terms:
\begin{equation}\label{eq:direct_adjoint_exp}
    \begin{pmatrix}
    \gamma_i^k\\
    \gamma_{i_1}^k
    \end{pmatrix}
     = \mathbb{E}^k\left[ D_i^k  \begin{pmatrix}
    \gamma_i^{k+1}\\
    \gamma_{i_1}^{k+1}
    \end{pmatrix}  \Big| {\cal V}^k \right]   + \frac{\rho}{N} \mathbb{E}^k\left[ \left( J_{i}^{k+1} + J_{i_1}^{k+1} \right) \begin{pmatrix}
    \partial \log h_i^k /\partial v_i^k\\
    \partial \log h_i^k /\partial v_{i_1}^k\
    \end{pmatrix}  \Big| {\cal V}^k \right]\,,
\end{equation}
where $\phi_i^M = \phi(v_i^M)$ with $\phi(v)$ given in the objective function~\eqref{eq:obj1}, and
\[
J_i^k := \mathbb{E}\left[ \phi_i^M | {\cal V}^k\right]\,,\quad \forall k = 0,1,\ldots, M-1\,.
\]
Next, we sample~\eqref{eq:direct_adjoint_exp}, as in the forward DSMC algorithm, over a single sample of the dynamics of $N$ discrete particles. We remark that $J_{i}^{k+1}$ and $J_{i_1}^{k+1}$ become $\phi_{i}^{M}$ and $\phi_{i_1}^{M}$ after sampling. Finally, we obtain the discrete adjoint equations, for any $i = 1,\ldots, N$, 
\begin{equation}\label{eq:direct_adjoint_sample}
    \begin{pmatrix}
    \gamma_i^k\\
    \gamma_{i_1}^k
    \end{pmatrix}
     =  D_i^k  \begin{pmatrix}
    \gamma_i^{k+1}\\
    \gamma_{i_1}^{k+1}
    \end{pmatrix}   +  \frac{\rho}{N} \left( \phi_{i}^{M} + \phi_{i_1}^{M} \right)  \frac{\partial \log h_i^k }{\partial v_i^k}\begin{pmatrix}
    1\\
    -1
    \end{pmatrix}\,.
\end{equation}
We used the fact that $ \frac{\partial }{\partial v_i^k} \log h_i^k  = - \frac{\partial }{\partial v_{i_1}^k} \log h_{i}^k $ based on the symmetry of elastic binary collision.  The ``final data" for $\gamma_i^k$ with $k=M$ is
\begin{equation}\label{eq:direct_adjoint_final}
\gamma_i^M =   \frac{\partial}{\partial v_i^M} \mathbb{E}[\bar{\phi}^M  |{\cal V}^M] =   \frac{\partial}{\partial v_i^M} \bar{\phi}^M = \frac{\rho}{N} \phi'(v_i^M)\,,
\end{equation}
but can change if the concrete form of the objective function varies.  In~\Cref{fig:DSMC}, we illustrate the adjoint relation between the forward DSMC  method (propagating forward in time) and the adjoint DSMC method (evolving backward in time) when we use the same randomness in the forward DSMC to sample~\eqref{eq:direct_adjoint_exp}.

Moreover, if the objective function involves the PDE solution at other time steps before the terminal time $T$, the adjoint equation~\eqref{eq:direct_adjoint_sample} needs to be modified. However, it can be derived in a similar way following the steps above. To sum up, the discrete adjoint particle system is solved backward in time starting with the final condition~\eqref{eq:direct_adjoint_final}, evolving based on~\eqref{eq:direct_adjoint_sample} till the initial time $t_0$. 

Finally, we take the derivative of the Lagrangian $L$~\eqref{eq:full_Lag} with respect to the parameter $\theta$, and obtain
\begin{equation*}
\partial_\theta L = \partial_\theta 
\left(  \sum_{i=1}^N \mathbb{E}_{\hat{v}_{0i} \sim  f_0(v;\theta)}  \left[ \gamma_i^0 \cdot \hat{v}_{0i} \right] \right) \,.
\end{equation*}
Since we restrict the Lagrangian over the solution manifold by performing the forward DSMC, the second and third terms in~\eqref{eq:full_Lag} disappear, and we have $\partial_\theta L = D_\theta J$, the true gradient of the objective function $J$; see~\Cref{subsec:adjoint} for more details.

To turn the gradient in its expectation form into a Monte Carlo gradient estimator, we have to examine how $\hat{v}_{0i}$ depends on the parameter $\theta$ in the initial condition $f_0(v;\theta)$. We often employ the \textbf{pathwise Monte Carlo gradient estimator} introduced in~\Cref{subsubsec:pw_grad}. That  is, we apply the ``reparameterization trick'':
\[
f_0(v;\theta) = T(v;\theta)\sharp g(v)\,,
\]
where the probability distribution $g(v)$ is parameter-independent and the push-forward map $T(v;\theta)$ inherits all the parameter dependence. As a result, we obtain the Monte Carlo gradient estimator
\[
D_\theta J \approx \sum_{i=1}^N
\gamma_i^0 \cdot  \partial_\theta T(\epsilon_i,\theta)\,,\quad \epsilon_i\sim g(v)\,,\quad  i = 1,\ldots,N\,.
\]

In this demonstration, we place the parameter dependence on the initial condition as an example, but in many setups, the parameter appears at other places of the Boltzmann equation, such as the collision kernel and the boundary condition. The adjoint equation system and the gradient formula will change, but again, they can be derived following the steps above, combined with the Monte Carlo estimation strategies introduced in~\Cref{sec:background}.

In summary, we have introduced the adjoint DSMC method to obtain the gradient of the Boltzmann equation-constrained optimization problem. In the derivation, we have employed the score function Monte Carlo gradient estimator (see~\Cref{subsubsec:score}) and the pathwise Monte Carlo gradient estimator (see~\Cref{subsubsec:pw_grad}). In the numerical validation of the adjoint DSMC method, the coupling method (see~\Cref{subsubsec:coupling}) can be used to provide a reference gradient value using the divided difference; see~\cite[Sec.~6]{caflisch2021adjoint}. All three primary Monte Carlo gradient methods are instrumental for this kinetic equation-constrained optimization problem.

\subsection{Numerical Examples}
In this section,   we show the results of numerical simulations computing the gradients of the objective function $ J(\theta) = \int_{\bbR^3} \phi(v) f(v, T) \rd v$, at the final time $t=T$ with respect to the parameter $\theta$ in the initial conditions $f_0(v;\theta)$.   The OTD approach yields an optimality condition before discretization. The continuous adjoint equation with respect to the Boltzmann equation~\eqref{eq:homoBoltz}  is 
\begin{equation}\label{eq:Boltz_adjoint}
  \left\{\begin{array}{rl}
- \partial_t  \gamma  &= \iint (  \gamma_1' + \gamma' - \gamma_1 - \gamma)   f(v_1,t)   q(v-v_1,\sigma)  \rd \sigma \rd v_1,\\
\gamma(v,T) &= - \phi(v). \\
\end{array}\right.
\end{equation} 
The gradient with respect to $\theta$ requires the solution to both the forward Boltzmann equation~\eqref{eq:homoBoltz} and the adjoint equation~\eqref{eq:Boltz_adjoint}, given as follows. 
\begin{equation}\label{eq:Bolz_OTD_grad}
\partial_\theta J = -\int\gamma(v,0) \partial_\theta f_0(v;\theta) \rd v\,.
\end{equation}

Four different methods for the gradient computation are used here: (i) finite difference method using several forward DSMC simulations with different parameter values, (ii) the adjoint DSMC method presented in Section~\ref{subsec:adjoint_DSMC},  (iii) the DSMC-type scheme for the continuous adjoint equation~\eqref{eq:Boltz_adjoint},  and (iv) the direct integration for solving the continuous adjoint equation~\eqref{eq:Boltz_adjoint}.  

Here we consider Maxwellian gas with a collision kernel $q(v-v_1,\sigma) =1/(4\pi)$ and assume $\rho(t)=\int_{\bbR^3} f(v,t) \rd v=1$.   Therefore,  $\mu =\rho \int_\Stwo q(\sigma)\rd \sigma = 1$.   We consider  $\phi(v)  = v_x^4$ where $v = [v_x, v_y, v_z]^\top$.   Regarding the parameter $\theta$,  we use temperature values in the initial distribution function $\theta=[T_x^0, T_y^0, T_z^0]$.   For all the methods,  we use the same initial condition,   an anisotropic Gaussian, 
\begin{equation*}~\label{eq:IC}
f_0(v) = \frac{1}{(2\pi)^{3/2}\sqrt{T_x^0 T_y^0 T_z^0}} \exp\left(-\frac{v_x^2}{2T_x^0}-\frac{v_y^2}{2T_y^0}-\frac{v_z^2}{2T_z^0}\right),
\end{equation*}
where $T_x^0=0.5,T_y^0=1,T_z^0=1$.  In all the tests, we use the forward Euler time-integration scheme with a time-step $\Delta t=0.1$.

\begin{table} %
\centering
{\renewcommand{\arraystretch}{1.2}%
\begin{tabular}{|c|c|c|c|c|}   \hline 
 \quad  &  $N = 10^6$  & $N = 10^7$ & $N = 10^8$ & $n_\text{grid}=30$\\
  \hline 
forward DSMC simulation &    $0.38$   sec & $5$ sec &  $60$ sec  & \\
adjoint DSMC simulation &  $0.22$ sec & $2.7$ sec & $30$ sec &  \\
the DSMC-type scheme &  $280$  sec  & $4100$ sec &\quad &  \\
direct integration of~\eqref{eq:Boltz_adjoint} &   &  & & $25000$ sec \\
   \hline
\end{tabular} }
  \caption{CPU run time of the four different methods approximating the gradient of a Boltzmann equation-constrained optimization problem~\cite{caflisch2021adjoint}.}\label{Table8}
\end{table}

\begin{figure}
\centering
\includegraphics[width=0.6\textwidth]{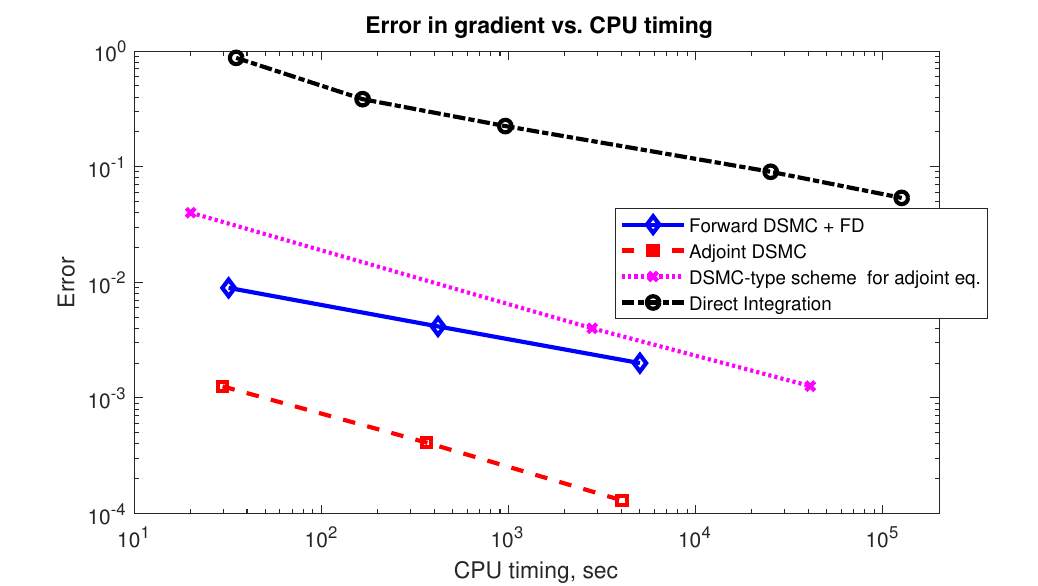}
\caption{The error in the gradient approximation vs. ~CPU time for all four methods~\cite{caflisch2021adjoint}.} \label{fig:Error_vs_CPU}
\end{figure}

All four methods lead to the same gradient values, but the adjoint DSMC method based on the DTO framework, as described in Section~\ref{subsec:adjoint_DSMC},  is the best in terms of performance given that we do not need many digits of accuracy.   Table~\ref{Table8} shows the CPU time of the code we recorded per objective function using Intel Core i7-3770K processor (4 cores @4.5Ghz) and $N=10^6,10^7,10^8$, $T=2$, $\Delta t=0.1$, $n_\text{grid}=30$ parameters. Here, $n_\text{grid}$ is the grid size per coordinate in the direct integration method (iv).   The adjoint DSMC is slightly faster (up to 20-30$\%$) than the forward DSMC for the same number of particles $N$.  The adjoint DSMC method (ii) is more than $10^3$ times faster than the DSMC-like scheme (iii) and much faster than the direct integration of~\eqref{eq:Boltz_adjoint},  which yields (iv).  At the same time, the errors in the adjoint DSMC are at least one order of magnitude smaller than in other methods due to the absence of finite-difference or interpolation errors.

Figure~\ref{fig:Error_vs_CPU} shows the numerical error in $\frac{\partial J}{\partial T_x^0}$ vs.  the~CPU time measured in simulations for all four methods~\cite{caflisch2021adjoint}.  Since one forward DSMC solve is needed for all methods, Figure~\ref{fig:Error_vs_CPU} reflects only CPU timing for the additional computations needed to compute the gradient, namely one extra forward DSMC simulation for (i) or one backward solution for (ii)-(iv).   Although the adjoint DSMC method is known to have a random error $\mathcal{O} (N^{-1/2})$,  it is much faster to evaluate than the DSMC-type scheme and the direct integration scheme to solve the continuous adjoint equation~\eqref{eq:Boltz_adjoint}.  This is an example where OTD-based approach is much less efficient than the DTO-based approach,  especially when the adjoint is computed at the level of particles.

\section{Conclusions}\label{sec:conclusion}
As described in the previous sections, Monte Carlo methods come with significant advantages, as well as difficulties: They are versatile and easy to implement, and they are not limited by dimensionality (more specifically, Monte Carlo methods typically depend linearly, not exponentially, on the dimension of an application); but they are slow and noisy (i.e., they have random errors). Indeed, their practical implementation often depends critically on noise control (i.e., variance reduction).

This review shows that the Monte Carlo method is useful for simulating adjoint equations or variables. This may be surprising because the differentiation that is inherent in adjoint methods could amplify the random noise in the Monte Carlo approach. Indeed, it is essential to avoid “differentiating the noise” in formulating adjoint Monte Carlo methods. This is done differently for OTD (optimize-then-discretize) or DTO (discretize-then-optimize) frameworks, in both of which random sampling is the discretization step (D), and differentiation occurs in the optimization step (O). For OTD, the differentiation occurs before sampling so that it does not involve random noise. On the other hand, for DTO, random sampling occurs first, and it is essential that each differentiation step only applies to the smooth influence of a single sample value and does not involve the differences between (random) sample values. 

A second theme of this review is that the remarkable effectiveness of the score function extends to adjoint Monte Carlo methods. In particular, we showed that an adjoint method can be applied to Monte Carlo sampling that involves a discontinuity due to the use of the rejection method, one of the most effective techniques for efficient sampling, in which the discontinuity is due to the choice between acceptance and rejection. We expect that similar use of the score function will be essential for further application of adjoint Monte Carlo methods.

The two application examples described in this survey --- radiative transport and the dynamics of rarefied gases and plasmas --- should be just the beginning. We expect that Monte Carlo adjoint methods will be broadly applicable to a wide range of models and simulation methods for physical, biological, and social sciences.

\section*{Acknowledgements}
We thank the editors of the special collection  ``Active Particles'' Vol 4., Jos\'e Carrillo and Eitan Tadmor, for their encouragement. R. Caflisch is partially supported by the Department of Energy (DOE) under grant DE-FG02-86ER53223. Y.~Yang is partially supported by the Office of Naval Research (ONR) under grant N00014-24-1-2088.

\bigskip
\bigskip

\end{document}